\documentclass[12pt]{article}
\begin{document}

\newtheorem{thm}{\noindent Theorem}[section]
\newtheorem{lem}{\noindent Lemma}[section]
\newtheorem{cor}{\noindent Corollary}[section]
\newtheorem{prop}{\noindent Proposition}[section]
\newtheorem{conj}{\noindent Conjecture}[section]
\newtheorem{assert}{\noindent Assertion}[section]

\makeatletter\renewcommand{\theequation}{%
\thesection.\arabic{equation}}
\@addtoreset{equation}{section}\makeatother

\newcommand{\qed}{\hbox{\rule[0pt]{3pt}{6pt}}}
\newcommand{\vlimsup}{\mathop{\overline{\lim}}}
\newcommand{\vliminf}{\mathop{\underline{\lim}}}
\newcommand{\Av}{\mathop{\mbox{Av}}}
\newcommand{\spec}{{\rm spec}}
\newcommand{\subarray}[2]{\stackrel{\scriptstyle #1}{#2}}
\setlength{\baselineskip}{14pt}

\def\textmc{\rm}
\def\({(\!(}
\def\){)\!)}
\def\R{{\bf R}}
\def\Z{{\bf Z}}
\def\N{{\bf N}}
\def\C{{\bf C}}
\def\T{{\bf T}}
\def\E{{\bf E}}
\def\H{{\bf H}}
\def\Prob{{\bf P}}
\def\M{{\cal M}}     
\def\F{{\cal F}}
\def\G{{\cal G}}
\def\D{{\cal D}}
\def\X{{\cal X}}
\def\A{{\cal A}}
\def\B{{\cal B}}
\def\L{{\cal L}}
\def\a{\alpha}
\def\b{\beta}
\def\e{\varepsilon}
\def\de{\delta}
\def\ga{\gamma}
\def\k{\kappa}
\def\la{\lambda}
\def\fa{\varphi}
\def\th{\theta}
\def\si{\sigma}
\def\t{\tau}
\def\om{\omega}
\def\De{\Delta}
\def\Ga{\Gamma}
\def\La{\Lambda}
\def\Om{\Omega}
\def\Th{\Theta}
\def\lan{\langle}
\def\ran{\rangle}
\def\lbr{\left(}
\def\rbr{\right)}
\def\const{\;\operatorname{const}} 
\def\dist{\operatorname{dist}} 
\def\Tr{\operatorname{Tr}}
\def\quadd{\qquad\qquad}
\def\n{\noindent}
\def\beq{\begin{eqnarray*}}
\def\eeq{\end{eqnarray*}}
\def\supp{\mbox{supp}}
\def\beqn{\begin{equation}}
\def\eeqn{\end{equation}}
\def\bp{{\bf p}}
\def\sg{{\rm sign\,}}
\def\br{{\bf r}}
\def\1{{\bf 1}}
\def\v2{\vskip2mm}
\def\n{\noindent}
\def\pf{{\it Proof.~}}

\begin{center}
{\Large  One dimensional lattice random walks with absorption \\  at a point /  on a half line}
\vskip6mm
{K\^ohei UCHIYAMA} \\
\vskip2mm
{Department of Mathematics, Tokyo Institute of Technology} \\
{Oh-okayama, Meguro Tokyo 152-8551\\
e-mail: \,uchiyama@math.titech.ac.jp}
\end{center}

\vskip6mm

\begin{abstract}
This paper concerns a random walk that moves  on the integer lattice and has zero mean and a  finite variance.
We obtain first an asymptotic estimate  of the transition probability  of    the walk  absorbed at the origin, and then, using the obtained estimate, that of the walk absorbed on a  half line. The latter is  used to evaluate the space-time distribution for the first entrance   of the walk  into the  half line. 
\footnote{
{\it key words}:  absorption, transition probability, asymptotic estimate, one dimensional   random walk}
 \end{abstract}
 
 \vskip6mm

\n
{\Large \bf Introduction}
\vskip2mm
Let $S^x_n=x+Y_1+\cdots+Y_n$ be a  random walk on the integer lattice $\Z$ starting at $x$ where the increments $Y_j$ are  independent and identically distributed  random variables  defined on some probability  space $(\Om, \F, P)$ and taking values in $\Z$.  Let $Y$ be a random variable having the same law as $Y_1$.  We suppose throughout the paper that the walk $S^x_n$ is irreducible and satisfies
\beqn\label{mom}
EY=0~~~~\mbox{and}~~~~\sigma^2:=E|Y|^{2}<\infty,
\eeqn
 where  $E$ indicates the expectation by $P$.
In this paper we compute an asymptotic form as $n\to\infty$ of the probability
\beqn\label{q0}
q^n(x,y)=P_x[S_n=y, S_1\neq 0,S_2\neq 0,\ldots, S_n\neq 0],
\eeqn
the transition probability of the walk absorbed at the origin, where (and in what follows) $P_x$ denotes the law of the walk $(S^x_n)_{n=0}^\infty$ and  under $P_x$ we   simply write $S_n$ for $S_n^x$.
The result on $q^n$ will be used  to  evaluate  $q_{(-\infty,0]}^n(x,y)$, the transition probability  
 of the walk that is absorbed  when it enters the negative half line, and  the result on the latter  in turn  to evaluate the space-time distribution for the first entrance   of $S^x_n$ into the negative half line.

The  local central limit theorem,  which gives  a  precise asymptotic form of  the transition probabilities $p^n(y-x):= P_x[S_n=y]$, plays a fundamental role in  both   theory and application  of  random walks,  
whereas concerning its analogue for $q_{(-\infty,0]}^n(x,y)$ or $q^n(x,y)$,   for all its significance,  there seem lacking,  except for   simple random walk case,   detailed results such as provide the precise asymptotic form  of them  [but see `Note added in proof'  at the end of the paper]. 
In this paper  we observe that the asymptotic forms of  both $q^n$ and $q_{(-\infty,0]}^n$  are  given by  the corresponding density of the  Brownian motion if  space variables $x, y$  as well as $n$ become large in a suitable  way, but  obviously they fails to be if $x$ and/or $y$ remain in a finite set. In the latter case  the order of magnitude of  decay (as $n\to\infty$)  does not differ   but the coefficients   do from  the  Brownian ones. These coefficients are expressed by means of  either the potential function of the walk  or  a pair of  \lq harmonic'  and \lq conjugate harmonic' functions on the positive half line  (renewal functions of ladder-height  processes) according as the absorption is made at the origin or on the negative half line.

 A primary estimate
 of $q^n$ is derived by using  Fourier analytic method; afterwards we   refine it by applying the result on the entrance distribution of $(-\infty,0]$ mentioned above (under an additional moment condition). Our results concerning $q_{(-\infty,0]}^n$  partly but significantly rest on a profound theory  of the random walk on the half line as found in Spitzer's book \cite{S}. 
The transition probability $q^n$ may be
viewed as the Green function of the space-time walk, an extremal case of two dimensional walks, absorbed on the coordinate axis of the time variable.  In a separate paper \cite{U3}  we study the corresponding problem for   two-dimensional random walks with zero mean and finite variances. With the help of some of the results obtained here and  in \cite{U3} the asymptotic estimates of the Green functions of 
the walks restricted on the upper half space are computed in \cite{U4}. 
A closely related issue  concerning  the  hitting distribution of a line for  two-dimensional walks is studied in \cite{Us}.

We illustrate  how fine  the estimate obtained is   by applying them to  a  problem on a system of independent random walks. Suppose that the  particles  are  initially placed on the positive half line of $\Z$ (one on each site) and   independently move according to the substochastic transition law $q^n$. Then  how does the total number of particles on the negative half line at time $n$ behave for large $n$?  We shall  prove that the expected number of such particles converges to a positive constant if $E[|Y|^3; Y<0]<\infty$ and diverges to infinity otherwise, provided that the walk is not left continuous; an analytic expression of the constant will be  given.

\section{Statements of Results}
 Let $S^x_n$ be the random walk described in Introduction and $P_x$ its probability law. Put $p^n(x)=P[S^0_n=x]$, $p(x)=p^1(x)$ and define the potential function 
 \beqn\label{a_def}
 a(x)=\sum_{n=0}^\infty[p^n(0)-p^n(-x)];
 \eeqn
 the series on the right side is convergent  and   $a(x)/|x|\to 1/\sigma^2$ as $|x|\to\infty$ (cf. Spitzer \cite{S}:Propositions P28.8
 and  P29.2). Denote by $d_\circ$ the period of the walk  (namely $d_\circ$ is the smallest positive integer such that $p^{d_\circ n}(0)>0$ for all sufficiently large $n$).
Put 
$${\sf g}_n(u)=\frac{e^{-u^2/2n_*}}{\sqrt {2\pi  n_*}}  ~~~~~~\mbox{where}~~~~~~n_*=\sigma^2 n.$$

The following notation will  also be used: $a\wedge b =\min\{a,b\}$, $a\vee b=\max \{a,b\}$ ($a,b\in \R)$; for functions $g$ and $G$ of a variable $\xi$,
 $g(\xi) = O(G(\xi))$  means  that there exists  a constant $C$ such that $|g(\xi)| \leq C|G(\xi)|$ whenever  $\xi$ ranges over a specified set; $ {\bf 1}({\cal S})$ denotes the indicator of a statement ${\cal S}$, i.e.,  
 $ {\bf 1}({\cal S})=1~\mbox{ or} ~~0~~\mbox{ according as ${\cal S}$ is true or not.}$
 
  We shall denote by $a_{\circ}$  an arbitrarily chosen  constant that is greater than unity (as in the items ${\bf (i)}$ and ${\bf (ii)}$ below), whereas positive constants to be determined  independently of variables $x, y, n$  etc.  but  may depend on the law of $Y$  are denoted by $C, C_1, C_2,\ldots$, whose values do not have to be the same at different occurrences  even though the same letter may be used.
 
 \v2\v2\n
{\bf 1.1.}~~ Let  $q^n(x,y)$ denote the transition probability of the walk  $S^x_n$ that is absorbed when it hits  the origin as defined by (\ref{q0}) 
(which entails that $q^n(x,y)=0$ if $y=0, x\neq 0$ and $q^0(x,y)={\bf 1}(x=y)$). For convenience sake we put 
$$a^*(x)=\1(x=0)+a(x).$$
 
\begin{thm}\label{thm1.1} ~
 The following asymptotic estimates of $q^n(x,y)$ as $n\to \infty$,  given in  three cases of constraints on $ x$ and $y$, hold true uniformly for $x$ and $y$ subject to the respective constraints. 

\v2
{\bf (i)}~  Under  $|x|\vee|y|< a_\circ \sqrt n$ and $|x|\wedge |y|=o(\sqrt n)$,
\beqn\label{(i)}
q^n(x,y)=\frac{\sigma^4a^*(x)a(-y)+xy}{n_*}\,p^n(y-x)+o\bigg(\frac{(|x|\vee1)|y|}{n^{3/2}}\bigg).
\eeqn

\v2
{\bf (ii)}~  Under  $a_\circ^{-1}\sqrt n < |x|,\,|y|< a_\circ \sqrt n$ (both $|x|$ and $|y|$ are between the two extremes),
\begin{eqnarray}\label{q(ii)}
q^n(x,y)&=&d_\circ {\bf 1}\Big(p^n(y-x)\neq 0)\Big)\Big[{\sf g}_{n}(y-x)-{\sf g}_{n}(y+x)\Big]+o\bigg(\frac{1}{\sqrt n}\bigg)~~~\mbox{if}~~xy>0,\\
q^n(x,y)&=&o\bigg(\frac{1}{\sqrt n}\bigg)~~~~~~~~~~~~~~~~~~~~~~~~~~~~~~~~~~~~~~~\mbox{if}~~~~~xy<0.
\label{(ii)}
\end{eqnarray}

\v2
{\bf (iii)}~   Let $0< |x|\wedge |y| <\sqrt n< |x|\vee|y|$.  ~ Then, if $E|Y|^{2+\de}<\infty$ for some $ \de\geq 0$, 
$$q^n(x,y)=O\bigg(\frac{|x|\wedge |y|}{|x|\vee|y|}{\sf g}_{4n}(|x|\vee|y|)\bigg)+o\bigg(\frac{|x|\wedge |y|}{(|x|\vee|y|)^{2+\de}}\bigg).
 $$
({\it ${\sf g}_{4n}$ on the right side can be replaced by ${\sf g}_{(1+\e)n}$ with any $\e>0$.})
\end{thm}

\v2\v2
As a simple consequence of {\bf (i)} and {\bf (iii)} of Theorem \ref{thm1.1} we have the bound
\beqn\label{iv}
q^{n}(x,y)\leq C\frac{(|x|+1)|y|}{n^{3/2}}.
\eeqn
valid  for all $ n, x$ and $y$.
It is  noted that if the walk is {\it left continuous}, i.e. $P[Y\le -2]=0$,  then $\sigma^2 a(x)=x$ for $x>0$, hence 
the leading term in the  formula of {\bf (i)} vanishes for $x>0, y<0$ in agreement with the trivial fact that $q^n$ itself does.  

If $E[|Y|^3]<\infty$ and $xy<0$, then  the assertion {\bf (i)}  can be  refined in two ways: the error term  in 
(\ref{(i)}) may be replaced by $o((|x|+|y|)n^{-3/2})$ and the resulting formula is valid  uniformly for $|x|\vee |y|<a_\circ \sqrt n$. Let $C^+$ be the constant given by
\beqn\label{C^+0}
C^{+}:= \lim_{x\to \infty}(\sigma^2 a(x)-x) \leq \infty.
\eeqn
We shall show (Corollary \ref{lem2.5} in Section 2; see also Corollary \ref{cor7.5})  that the limit exists
and that it  is finite if and only if  $E[|Y|^3; Y<0]<\infty$ and positive unless the walk is left continuous.  It follows that 
$$\sigma^4a^*(x)a(-y)+xy= C^+(x-y)(1+o(1))~~~\mbox{ as}~~ x\wedge (-y)\to \infty,$$
 provided $E[|Y|^3; Y<0]<\infty$. In view of this  relation and duality the refined version of   {\bf (i)} mentioned above  may  read as follows.
\begin{thm}\label{thm1.5}~
Suppose that   $E[|Y|^3;Y<0]<\infty$. Let $y<0<x$. Then uniformly for $x\vee|y|\le a_\circ\sqrt n$, as  $x\wedge |y|\to \infty$
\beqn
q^n(x,y)=C^+\frac{x+|y|}{n_*}p^n(y-x)+ o\bigg(\frac{x\vee|y|}{n^{3/2}}\bigg).
\eeqn
\end{thm}

\v2
The proof of Theorem \ref{thm1.5}  requires more delicate analysis than that of  Theorem \ref{thm1.1};  it rests  on Theorem \ref{thm1.4} below and will be given  after the proof of it 

Given a constant  $\a\in (0,1)$,  one may consider the  absorption  which is not absolute but takes place with  probability $\a$ each time the walk is about to visit the origin.  Denote by $q_\a^n(x,y)$ the transition probability of the  process
subject to such absorption.   In  Section 6 we shall obtain the  asymptotic estimates  
\beqn\label{qqq}
q_\a^n(x,y)-q^n(x,y)=\frac{(1-\a)\sigma^2}{\a}\cdot\frac{a^*(x)+a^*(-y)}{n}p^n(y-x)(1+o(1))
\eeqn
valid uniformly for $|x|\vee|y|<a_\circ \sqrt n$.
Note that  as $|x|\wedge |y|\to\infty$ under the same constraint on $x, y$,  the right side divided by $q^n(x,y)$ tends to zero  for $xy>0$, while it is asymptotically a positive constant  for  $xy<0$,  provided $E[|Y|^3]<\infty$ according to  Theorem \ref{thm1.5}.

\v2\v2\n
{\bf 1.2.}~~ Here we consider the walk  absorbed when it enters $(-\infty,0]$.  Let $T$ denote  the first entrance time into $(-\infty,0]$:
 $$T=\inf\{n\geq 1: S_n\leq 0\},$$
 and   $q_{(-\infty,0]}^n(x,y)$ the transition probability of the absorbed walk:  
$$q_{(-\infty,0]}^n(x,y)=P_x[S_n=y, S_1>0,\ldots, S_n > 0]= P_x[S_n=y, n<T]~~~~ ~~(x,y >0).$$

The next result states that  $q_{(-\infty,0]}^n(x,y)$  behaves  similarly to $q^n(x,y)$ within any parabolic region if both  $x$ and $y$ get large. 

\begin{prop}\label{thm1.2}~ Uniformly for $n\geq (x\vee y)^2/a_\circ$, as $x\wedge y\to \infty$
$$q_{(-\infty,0]}^n(x,y)=q^n(x,y)(1+o(1)).$$
\end{prop}

Let $f_+(x)$ (resp. $f_-(x)$) ($x= 1,2,\ldots$) be the positive  function on $x> 0$  that is asymptotic to $x$  as $x\to\infty$ and harmonic with respect to  the walk  $S_n$ (resp $-S_n$) absorbed on $(-\infty, 0]$ : 
\beqn\label{f_def}
~~f_{\pm}(x)=E[f_{\pm}(x\pm Y);\, x\pm Y >0]~~(x\geq 1)~~~\mbox{and}~~~ \lim_{x\to\infty} f_{\pm}(x)/x=1,
\eeqn
each of  which exists uniquely (Spitzer \cite{S}:P19.5).  (It is warned  that it is not  $[1,\infty)$  but   $[0,\infty)$  on which  the harmonic function is considered in \cite{S}.) 
 
\begin{thm}\label{thm1.3}~ Uniformly for $0< x, y\leq a_\circ\sqrt n$, as $xy/n\to 0$
$$q_{(-\infty,0]}^n(x,y)=\frac{2f_+(x)f_-(y)}{n_*}p^n(y-x)(1+o(1)).$$
\end{thm}

 \vskip2mm
 From Theorem \ref{thm1.3} one derives an asymptotic form of the space-time  distribution of  the first entrance into $(-\infty,0]$, which we denote by $h_x(n,y)$:  for $y\leq 0$
 $$h_x(n,y) =P_x[S_T=y, T=n].$$
Put 
\beqn\label{q}
H_{\infty}^+(y)=\frac2{\sigma^2}E[f_-(y-Y);Y<y]=\frac{2}{\sigma^{2}}\sum_{j=1}^\infty f_-(j)p(y-j)~~~~~ (y\leq 0).
\eeqn

\begin{thm}\label{thm1.4}~ Suppose $E[|Y|^{2+\de};Y<0]<\infty$ for some $\de\ge 0$ and $d_\circ=1$. Then, uniformly for $y\leq 0< x\leq a_\circ\sqrt n$, as $n\to \infty$
\beqn\label{h0}
h_x(n,y)=\frac{f_+(x){\sf g}_n(x)}{n}H_{\infty}^+(y)(1+o(1))+\frac{x}{n^{3/2}}\a_n(x,y),
\eeqn
with
$$\a_n(x,y)=o\Big((|y|\vee\sqrt n\,)^{-1-\de}\,\Big),~~~\sum_{y\leq 0}|\a_n(x,y)|=o(n^{-\de/2}) ~~~~\mbox{and}~~~\sum_{y\leq 0}|\a_n(x,y)||y|^\de=o(1);$$
and for   $x\geq \sqrt n$ and $y\leq 0$
\beqn\label{upb-h20}
h_x(n,y)\leq C\bigg[\frac{{\sf g}_{4n}(x)}{\sqrt n}+o\bigg(\frac1{x^{2+\de}}\bigg)\bigg]H_{\infty}^+(y) +\frac{C}{\sqrt n}P[Y<y-{\textstyle \frac12}x],
\eeqn
and in particular
\beqn\label{eq1.4}
h_x(n,y)\leq C{H_\infty^+(y)}x^{-1} n^{-1/2}.
\eeqn
\end{thm}

 Since $P_x[\,T=n]=\sum_{y\leq 0}h_x(n,y)$ we have the following corollary of Theorem \ref{thm1.4}.
\begin{cor}\label{cor1.1}~Uniformly in $x\geq 1$
$$P_x[\,T=n]= \frac{f_+(x){\sf g}_n(x)}{n}(1+o(1)) + o\bigg(\frac{x}{n^{3/2}}\wedge \frac1{x\sqrt n}\bigg).$$
\end{cor}

\vskip2mm  
{\sc Remark.}~  (a)~ $H_{\infty}^+$ is the probability on $(-\infty,0]$  that arises as  the limit as $x\to\infty$ of the first entrance distribution $H^+_x(\cdot)=\sum_n h_x(n,\cdot)$ (\cite{S}, P19.4). 
This  in particular gives the identity
$\sum_{j=1}^\infty f_-(j)P[Y\leq -j]=\sigma^2/2.$
 
(b) ~ If the walk is right continuous (i.e., $P[Y\geq 2]=0$) as well as in the case  when it is left continuous we have $q^n_{(-\infty,0]}(x,y)=q^n(x,y)$ for $ x, y>0$.

(c) ~   The formula (\ref{h0})  holds true also in the periodic case  (i.e., $d_\circ>1$),  if    the leading term on its right side is  multiplied by $d_\circ {\bf 1}\Big(p^n(y-x)\neq 0)\Big)$ as in (\ref{q(ii)}).

(d)~ The function $f_-$ may be given by the formula
$$f_-(x)=f_-(1)\Big (1+ E_0[\,\mbox{the number of ascending ladder points} \in [1, x-1] \,]\Big)$$
and its dual formula  for $f_+$ (\cite{S}:pp.201-203). Under   our normalization of $f_{\pm}$ the initial value $f_-(1)$ (resp. $f_+(1)$) equals the expectation of the strictly ascending (resp. descending)  ladder height : 
\beqn\label{f(1)}
f_-(1)=E_0[S_{\tau([1,\infty))}]~~~ \mbox{and}~~~ f_+(1)=-E_0[S_{\tau((-\infty,-1])}],
\eeqn 
where $\tau(B)$ denotes the first entrance time into a set $B$, in view of the renewal theorem.

(e) ~If the starting point is 1,  the  Baxter-Spitzer  identity  gives
\beqn\label{h_1}
\sum_{n=0}^\infty r^n\sum_{y\leq 0} z^{1-y}h_1(n,y)=1-\exp \bigg(-\sum_{k=1}^\infty \frac{r^k}{k} E_0[z^{-S_k}; S_k < 0]\bigg) ~~~~~~(|z|\leq 1, |r|<1)
\eeqn
 and a similar formula for $q_{(-\infty,0]}^k(1,y) $  (\cite{C}:Theorem 8.4.2, \cite{F}: Lemmas 1 and 2 of Section XVIII.3, \cite{S}:P17.5). We shall  use these identities  not directly but via  certain  fundamental results (including those on $f_{\pm}$ and found in \cite{S}) that  are based on them. Taking $z=1$ 
 the above formula reduces to
 $$1-E_1[r^T]=\sqrt{1-r} \exp \bigg[\sum_{k=1}^\infty \frac{r^k}{k} \bigg(\frac12-P_0[S_k<0] \bigg)\bigg],$$
  and, applying  Karamata's Tauberian theorem,  one can readily find an asymptotic formula of $P_1[T\geq n]$, which is also obtained from Corollary   \ref{cor1.1} and (\ref{f(1)}). It would however  be  difficult to derive directly  from   the formula (\ref{h_1})  such  fine   estimates of $h_1(n,y)$ as given in Theorem \ref{thm1.4}.

\v2\v2\n
{\bf 1.3.}~~ 
 For $x\in \Z$, let $Q^+_x(n)$ denote the probability that the walk starting at $x$  is found in the negative half line at time $n$ without having hit the origin  before $n$:
$$Q_x^+(n)=\sum_{y=-\infty}^{-1}q^n(x,y).$$
\begin{prop} \label{prop1.3.1}~~As $x/\sqrt n\to 0$
\beqn\label{Q+}
Q_{x}^+(n)=\frac{\sigma^2a^*(x)-x}{\sqrt{2\pi n_*}}+o\bigg(\frac{|x|+1}{\sqrt n}\bigg);
\eeqn
and uniformly in $n$, as $x\to\infty$
\beqn\label{Q+-}
Q_{-x}^+(n)=\int_{-x}^x {\sf g}_n(u)du\Big[1+o(1)\Big].
\eeqn
If $E[|Y|^{3}; Y<0]<\infty$, then for $x>0$,  the error term in (\ref{Q+}) can be replaced by $o(1/\sqrt n)$. 
\end{prop}
 
 The formula (\ref{Q+-}) follows from (\ref{Q+}) if $x/\sqrt n\to 0$ so that it signifies only in the case $x>a_\circ^{-1}\sqrt n$.

Let $C^+$ be the same constant as introduced in the subsection 1.1 (just before Theorem \ref{thm1.5}). $C^+$ is finite if and only if $E[|Y|^3;Y<0]<\infty$ as remarked there.
\v2 
\begin{thm} \label{thm1.3.2}~  Let $ \nu_n= \sum_{x=1}^\infty Q^+_{x}(n)$.  Then
${\displaystyle \lim_{n\to\infty} } \nu_n=\frac12 C^+.$
\end{thm}

 One can extend Theorem \ref{thm1.3.2} as follows.  We are concerned with the particles each of which performs random walk according to the transition law $q^n(x,y)$  independently of the other ones. Consider an experiment such that  at a time $n$ that is determined  prior to the experiment  the experimenter counts the number of particles lying in any interval of the negative half line for the system of our particles in which  at time 0 the particles are randomly placed at each site $x>0$ whose mean number,  denoted by $m_n(x)$, may depend on $n$ as well as $x$.  Let $N_n(\ell)$, $\ell>0,$ denote the number of particles that are found in the interval $[-\ell \sqrt {n_*},-1]$ at the  time $n$. The following extension is a corollary of the proofs of Proposition \ref{prop1.3.1} and Theorem \ref{thm1.3.2}.
 
 \begin{cor}\label{cor1.2} ~
 Suppose $m_n(x)=0$ for $x<0$, $m_n(x)$ is uniformly bounded and for each $K>0$, $m_n(x)=1+o(1)$  as $ n\to\infty$ uniformly for $0<x<K\sqrt n$. Then for  each positive  number $\ell $,
\beqn\label{experiment}
\lim_{n\to\infty} E[N_n(\ell)]=\frac{C^+}{\sqrt{2\pi}}\int_0^\ell e^{-t^2/2}dt.
\eeqn
\end{cor}

\v2\v2
 That  $\nu_n= \sum_{x=1}^\infty Q^+_{x}(n)$ is bounded if and only if  $E[|Y|^3; Y<0]<\infty$ is easy  to prove. Indeed 
  $$\nu_n = \sum_{k=1}^n \sum_{w=1}^\infty\sum_{y=-\infty}^{-1}Q^{*+}_w(k-1)p(y-w)Q^+_y(n-k),$$
 where $Q^{*+}_w(k)=\sum_{x=1}^\infty q_{(-\infty,0]}^k(x,w)$.  Crude
 applications of  Theorem \ref{thm1.1} and Proposition  \ref{thm1.2}  give
 \beqn\label{Q^-}
 C_1{\bf 1}\bigg(-1<\frac{y}{\sqrt n} <0\bigg)\frac{|y|}{\sqrt n} \leq Q^+_y(n) \leq C_2\frac{|y|}{\sqrt n}~~~~~\mbox{for}~~~y<0
 \eeqn
 and similar bounds for $Q^{*+}_u(n)$, respectively, whereupon,
noting  $\sum (k(n-k))^{-1/2} \sim \int_0^1 (t(1-t))^{-1/2}dt$, one finds that $\nu_n$ is bounded if and only if
 $\sum_{w=-\infty}^{-1} \sum_{y=-\infty}^{-1} p(y+w) wy<\infty$, but the latter condition  is equivalent to
 $E[|Y|^3; Y<0]<\infty$. 

\v2
The rest of the paper is organized as follows. In Section 2 we give some preliminary lemmas. The estimation of $q^n$ and that of $q^n_{(-\infty,0]}$ and $h_x(n,y)$  are carried out in Section 3 and Section 4, respectively.  Further detailed estimation of $q^n(x,y)$ for $xy<0$ that leads to the proof of Theorem \ref{thm1.5} is made in the end of Section 4. In Section 5 $Q^+_x(n)$ is dealt with. In Section 6 we briefly discuss on
$q^n_\a(x,y)$ and prove (\ref{qqq}).

\section{ Preliminary Lemmas }
This section  is divided into four  subsections.  In the  first one we give  some terminologies and notation as well as  some  fundamental results from Spizer's book \cite{S} in addition to those given in Section 1.
Both the  the second and the third ones depend  in an essential way on the classical results given in  the first subsection   but self-contained otherwise.
 \vskip4mm\n
{\bf 2.0.}~
 Let   $B$  be a subset of $\Z$. 
Denote by $\tau_B$  the first time when  $S_n$ enters $B$  after time $0$; $\tau_B=\inf\{n\geq 1: S_n\in B\}$.  For a point  $x\in \Z$ write $\tau_x$ for $\tau_{\{x\}}$.  For typographical reason we sometimes write
$\tau(B)$ for $\tau_B$.   

A function $\fa(x)$ on $\Z\setminus B$ that is bounded from below is said to be {\it harmonic} on $\Z\setminus B$ if $E_x[\fa(S_1): S_1\notin B]=\fa(x)$ for all $x \notin B$.  From this property 
with   the help of   Fatou's lemma   one infers  that for any Markov time $\tau$, $E_x[\fa(S_\tau); \tau<\tau_B]\leq \fa(x)$~ ($x \notin B$).  The functions $f_+(x)$  and $a(x)$ introduced in Section 1 (see (\ref{f_def}) and (\ref{a_def})) are harmonic on $[1,\infty)$ and on $\Z\setminus \{0\}$, respectively (\cite{S}, T29.1). The function $f_-(x)$ is harmonic also on $[1,\infty)$ but for the {\it dual} walk, namely the walk determined by the probability  law $p^*(x)=p(-x)$.

Let  $g_{\,(-\infty,0]}(x,y)$ ($\,x,y>0\,$)  denote the Green function of the walk $S_n$  absorbed on $(-\infty,0]$: $g_{\,(-\infty,0]}(x,y)=\sum_{n=1}^\infty g^n_{(-\infty,0]}(x,y)=\sum_{n=0}^\infty P_x[S_n=y, n<T]$, where $T=\tau_{(-\infty,0]}$ as in Section 1.  
It follows from the propositions P18.8, P19.3, P19.5 of \cite{S} that  the  increments
 $$u^{\pm}(y):=f_{\pm}(y)-f_{\pm}(y-1)~~~ (y=1,  2,  \ldots), ~u^\pm(0):=0$$
are all positive and have limits $\lim_{y\to\infty}u^{\pm}(y)= 1$ and with them the function $g_{\,(-\infty,0]}$ is expressed as
\beqn\label{g2}
g_{\,(-\infty,0]}(x,y)= \frac{2}{\sigma^2}\sum_{z=0}^{x\wedge y} u^+(x-z)u^-(y-z)~~~~~ (x,y >0).
\eeqn

Similarly  let $g_{\{0\}}(x,y)$ be the Green function of the process $S_n$ absorbed at the origin: $g_{\{0\}}(x,y)=\sum_{n=0}^\infty q^n(x,y)$.  Then, according to Spitzer \cite{S} (P29.4)
\beqn\label{g}
g_{\{0\}}(x,y)=a(x)+a(-y)-a(x-y)~~~~~~(x, y \in \Z\setminus \{0\}).
\eeqn
The results given in the following subsections ${\bf 2.1}$ and ${\bf 2.2}$, though easy consequences
of (\ref{g2}) and (\ref{g}),  do not seem to appear in the existing literature.

\vskip2mm\n
{\bf 2.1.}~ Let $H_x^+(y)$ denote (as in {\sc Remark} (a)) the  hitting distribution of $(-\infty,0]$ for the walk $S_n^x$:
\beqn\label{hd}
H_x^+(y):=P_x[S_{T}=y]~~~~~ (x>0, y\leq 0),
\eeqn 
which  may be  expressed as
\beqn\label{h^-}
 H_x^+(y)=\sum_{w=1}^\infty g_{(-\infty,0]}(x,w)p(y-w).
\eeqn
In view of  (\ref{g2}) we have $g_{(-\infty,0]}(x,w)\leq C f_-(w)$, hence
\beqn\label{17}
H_x^+(y)\leq CH_{\infty}^+(y).
\eeqn

\begin{lem} \label{lem2.50} ~~For $x>0$ and $y\leq 0$,
$${\rm (a)}~~~~~~~~~~~~~~~~~~~~~~~~~~~~~~~\sum_{z=-\infty}^{0} H_x^+(z)a(z-y)=a(x-y)-\sigma^{-2}f_+(x).~~~~~~~~~~~~~~~~~~~~~~~~~~~~~~~~$$
$${\rm (b)}~~~~~~~~~~~~~~~~~~~~~~~~~~~~~~~~~\sum_{z=-\infty}^{0} H_{\infty}^+(z)a(z)=\lim_{x\to\infty}\Big[a(x)-\sigma^{-2}f_+(x)\Big].~~~~~~~~~~~~~~~~~~~~~~~~~~~~~~~~~~$$
(Both sides of {\rm (b)} may be infinite simultaneously.) 
\end{lem}
\v2\n
\pf  ~ With $y\leq 0$ fixed define $\fa(x)=\sum_{z=-\infty}^{0} H_x^+(z)a(z-y)$ for $x>0$ and $\fa(x)=a(x-y)$ for $x\le 0$.  Owing to (\ref{g2}) and (\ref{h^-}) 
\beqn\label{eq2.5}
H_x^+(z)\leq C\sum_{w=1}^\infty (x\wedge w)p(z-w)\leq CxP[Y<z]~~~~~~~(y<0< x),
\eeqn
which combined with  $\sum_{z\leq 0} |z|P[Y<z]\leq\sigma^2<\infty$  shows that $\fa(x)$ takes a finite value; moreover, by  dominated convergence, $\fa(x)/x\to 0$ as $x\to \infty$.  It is  observed that $\sum_{z=-\infty}^{\infty}p(z-x)\fa(z)=\fa(x)$  for $x>0$ and $\sum_{z=-\infty}^{\infty}p(z-x)a(z-y)=a(x-y)$  for $x\neq y$.  Hence  $a(x-y)-\fa(x)$, vanishing  on  $x\leq 0$,  is  harmonic on  $x>0$ and asymptotic to $x/\sigma^2$ as $x\to\infty$. 
We may now conclude that   $a(x-y)-\fa(x)$  agrees with $\sigma^{-2}f_+(x)$ for $x>0$ since  a harmonic function on $x>0$ that is bounded below is unique apart from a constant factor.  Thus (a) has been verified. 
Let  $y=0$ and $x\to\infty$ in (a).  If the left side of (b) is infinite, so is the right side in view of  Fatou's lemma.  If it is finite,  the dominated convergence theorem  may apply owing to  (\ref{17}). ~~~\qed
\v2

 The proof of Lemma \ref{lem2.50} may be repeated word for word  but with $a(z-y)$ replaced by $z$ to yield 
\beqn\label{eq2.51}
\sum_{z=-\infty}^{-1} H_x^+(z)|z|= f_+(x)-x
\eeqn
and
\beqn\label{eq2.50}
\sum_{z=-\infty}^{-1} H_{\infty}^+(z)|z|=\lim_{x\to\infty}\Big[f_+(x)-x\Big].
\eeqn
Taking $y=0$ in  (a)  of Lemma \ref{lem2.50}  and combining it with (\ref{eq2.51}) we obtain 
\beqn\label{22}
\sum_{z=-\infty}^{-1} H_x^+(z)(\sigma^2 a(z)-z)=\sigma^2a(x)-x~~~~~~(x>0).
\eeqn

 Here we advance a corollary of Lemma \ref{lem2.50} that involves the constant  $C^+$ introduced in the subsection 1.1. It is convenient to define it by
$$C^+=\sum_{y=-\infty}^{-1}H_{\infty}^+(y)\Big[\sigma^2a(y)+|y|\Big]$$
rather than by (\ref{C^+0}). The relation (\ref{C^+0}) then ensues  as stated in the  corollary below. According to this definition it is clear that
 $C^+$ is finite if and only if $E[|Y|^3;Y<0]<\infty$,  and  positive unless the walk is left continuous. 
Define $H_{-\infty}^-(y)$ and $C^-$ analogously to $H_{\infty}^+$ and $C^+$:
$$H_{-\infty}^-(y)=\frac2{\sigma^2}E[f_+(Y-y);Y>y]~~~ (y\geq 0)~~~\mbox{and}~~~C^-=\sum_{y=1}^\infty H_{-\infty}^-(y)(\sigma^2 a(y)+y).$$
It holds that $C^-<\infty$ if and only if $E[|Y|^3; Y>0]<\infty$.

\begin{cor} \label{lem2.5} ~
${\displaystyle C^+=\lim_{x\to+\infty}(\sigma^2 a(x)-x)~\mbox{~and~}~C^-=\lim_{x\to -\infty}(\sigma^2 a(x)-|x|).}$~
\end{cor}
\vskip2mm\n
{\it Proof}.~  From (\ref{eq2.50}) and
 (b) of Lemma \ref{lem2.50} one deduces  the first relation  of the corollary.  The second one is   its dual. 
 ~~~\qed

\vskip4mm\n
{\bf 2.2.}~ The results in this subsection  are  somewhat different in nature from and independent of those of the preceding one (except for the use of (\ref{17}) ) although  machinery for the proof is essentially the same. Recall that  $T$ is written for $\tau_{(-\infty,0]}$. We shall show that $P_x[\tau_{[N,\infty)}<\tau_{0}]$ and  $P_x[\tau_N<\tau_0]$ are asymptotically equivalent (see Proposition \ref{lem2.2}). For the moment we obtain the following
\begin{lem}\label{lem2.1}~~Uniformly in $0<x<N$, as $N\to\infty$
$$\frac{a(x)}{a(N)}\geq P_x[\tau_{[N,\infty)}<\tau_{0}]\geq P_x[\tau_N<\tau_0]=\frac{\sigma^2a(x)+x}{2N}(1+o(1)).$$
\end{lem}
\v2\n
{\it Proof.}   We have $P_x[\tau_N<\tau_0]=g_{\{0\}}(x,N)/g_{\{0\}}(N,N)$ and the last relation of the lemma follows from (\ref{g}) together with
$$g_{\{0\}}(N,N)=a(N)+a(-N)=\frac{2N}{\sigma^2}(1+o(1)),~~~~\lim_{y\to\infty}[a(-y)-a(1-y)]=1/\sigma^2$$
(\cite{S}: P29.2). Since $a(x)$ is positive and harmonic on $\Z\setminus\{0\}$ (i.e. $E_x[a(Y+x)]=a(x),~x\neq 0$) and non-decreasing for $x>0$ large enough, 
$$a(x)\geq  E_x[a(S_{\tau([N,\infty))});\tau_{[N,\infty)}< \tau_0]\geq {a(N)}P_x[\tau_{\,[N,\infty)}< \tau_0].$$
Thus $P_x[\tau_{\,[N,\infty)}<\tau_0]\leq a(x)/a(N)$ provided $N$ is large enough, verifying the first inequality of the lemma. The second one is trivial. ~~~\qed 
\begin{lem}\label{unif_bd}~
As $N\to\infty$
\[\label{x}  \sup_{z>N}P_{z}[ \tau_N> T] \asymp \bigg[ N^{-1}\sum_{y=1}^{N-1}yH^+_{\infty}(-y)+\sum_{y=N}^\infty H^+_{\infty}(-y)\bigg] \longrightarrow 0.
\]
\end{lem}
\vskip2mm
{\it Proof.}~
We use the decomposition
$$
P_z[ \tau_N> T]=\sum_{w<N}P_z[S_{\tau((-\infty,N] )}=w]\Big(P_w[T<\tau_N]{\bf 1}(w>0)+{\bf 1}(w\leq 0)\Big).
$$
Writing the first probability under the summation sign by means of $H_x^+(y)$ (defined in (\ref{hd}))
and using the bound
 $H_x^+(y) \leq CH_\infty^+(y)$
 (see (\ref{17}))  together with Lemma \ref{lem2.1}  we obtain
 $$P_z[ \tau_N> T]  \leq\frac{C}{a(-N)}\sum_{-N<y<0}a(y)H_{\infty}^+(y) +\sum_{y\leq -N}H_{\infty}^+(y),$$
the right side  approaching zero. The lower bound is obtained by an application of Fatou's lemma.~~\qed
\begin{prop}\label{lem2.2}~~Uniformly in $0<x<N$, as $N\to\infty$
\vskip3mm\n
{\rm (a)}~$\displaystyle~~~~~~~~~~~~~~~~~~~~~~~~ P_x[\tau_{\,[N,\infty)}<T\,] - P_x[\tau_N<T\,]=o(x/N);~~~~~~~~~~~~~~~~~~~~~~~~~~~~~~~~~~~~~~~~~~~~~~$
$${\rm (b)}~~~~~~~~~~~~~~~~~~~~~~~~~ P_x[\tau_{\,[N,\infty)}<\tau_0\,] - P_x[\tau_N<\tau_0\,]=o(x/N).~~~~~~~~~~~~~~~~~~~~~~~~~~~~~~~~~~~~~~~~~~~~~~$$
\end{prop}
\vskip2mm
{\it Proof.}~~ The difference on the left side of (a) is  expressed as
\beqn\label{eq-1}
\sum_{z>N} P_x[\tau_{\,[N,\infty)}< T,\, S_{\tau([N,\infty))} =z]P_z[ \tau_N> T],
\eeqn
and hence  (a)   follows from the  preceding two  lemmas  (and the inequality $T\leq \tau_0$).  

The proof of (b) is similar: one has only to replace $T$  by $\tau_0$ in (\ref{eq-1}).  ~~~\qed
\vskip2mm

The next Proposition refines Theorem 22.1 of \cite{S} where the problem is treated by a quite different method from the present one.
\begin{prop}\label{lem2.3}~~Uniformly for $1\leq x<N$, as $N\to \infty$  
$$ P_x[ \tau_{\,[N,\infty)}<T\,] =  \frac{f_+(x)}{N} +o\bigg(\frac{x}{N}\bigg).
$$
\end{prop}
\vskip2mm  
\pf~~
By Proposition \ref{lem2.2} 
$$P_x[\tau_{\,[N,\infty)}<T\,]=P_x[\tau_N<T\,]+o\bigg(\frac{x}{N}\bigg)=\frac{g_{\,(-\infty,0]}(x,N)}{g_{\,(-\infty,0]}(N,N)}+o\bigg(\frac{x}{N}\bigg).$$
It is readily inferred that  as $N\to \infty$  
$$g_{\,(-\infty,0]}(x,N)=f_+(x)\bigg(\frac{2}{\sigma^2} +o(1)\bigg)  ~~~\mbox{uniformly for}~ ~1\leq x\leq N,$$
in particular, $~ g_{\,(-\infty,0]}(N,N)=2\sigma^{-2}N+o(N)$ and substitution leads to the desired relation. ~~~\qed

\begin{prop}\label{lem2.4}~~Uniformly for $1\leq x<N$, as $N\to \infty$  
$$\frac1{x}E_x[\,S_{\tau([N,\infty))};\, \tau_{\,[N,\infty)}<T\,]=\frac{N}{x}P_x[ \tau_{\,[N,\infty)}<T\,] +  o(1),
$$
or, what is the same thing, ~$E_x[\,S_{\tau([N,\infty))}-N|\, \tau_{\,[N,\infty)}<T\,]=o(N)$.
\end{prop}
\vskip2mm
{\it Proof.}~~That $f_+$ is non-negative and harmonic on $[1,\infty)$ implies that for $x>0$,
$$E_x[f_+(S_{\tau([N,\infty))}); \tau_{\,[N,\infty)}<T\,]\leq f_+(x).$$
Hence, employing  Lemma \ref{lem2.3},  one first observes that  uniformly for $1\leq x<N$
\beq
0\leq E_x[\,f_+(S_{\tau([N,\infty))})-f_+(N)\,; \tau_{\,[N,\infty)}<T\,] &\leq& f_+(x)-f_+(N)P_x[\tau_{\,[N,\infty)}<T\,]\\
&=& f_+(x)(1- f_+(N)/N)+o(x)
\eeq
 and then use  $\lim f_+(x)/x=1$ to find the formula of the proposition.~~~\qed

\vskip2mm
\begin{lem}\label{lem2.9}~  Uniformly for $0<x<N$, as $N\to\infty$
\beqn\label{23}
P_x[\tau_{ [N,\infty)}<\tau_0] =\frac{\sigma^2 a(x)+x}{2N}(1+o(1));~~\mbox{and}
\eeqn
$$\frac{\sigma^2 [a(x)+a(N)-a(x+N)]}{2N}(1+o(1))\leq P_x[\tau_{ (-\infty,-N]}<\tau_0] \leq\frac{\sigma^2 a(x)-x}{2N}(1+o(1)).~~~~~~
$$
\end{lem}
\vskip2mm
{\it Proof.}~ The first relation (\ref{23}) follows from (b) of Proposition \ref{lem2.2} and the last equality in Lemma \ref{lem2.1}.  The lower bound of the second relation is obtained in the same way as  is the second inequality in Lemma \ref{lem2.1}.  For the upper bound we apply (\ref{23}) (or rather its dual) and (\ref{22}) in turn to see 
\beq
P_x[\tau_{ (-\infty,-N]}<\tau_0] &=& \sum_{y= -N+1}^{-1}H_x^+(y)P_y[\tau_{(-\infty,-N]}<\tau_0]+\sum_{y\leq-N} H_x^+(y)\\
~~~~~&&\leq \sum_{y=-\infty}^{-1} H_x^+(y)\frac{\sigma^2 a(y)-y}{2N}(1+o(1)) = \frac{\sigma^2 a(x)-x}{2N}(1+o(1)).
\eeq
The proof of the lemma is complete.  ~~ \qed
\begin{prop}\label{prop2.4}~~ Uniformly for $0< |x|<N$, as $N\to\infty$
$$P_x[ \tau_{\Z\setminus (-N, N)} <\tau_0]= \frac{\sigma^2 a(x)}{N}(1+o(1)).
$$
\end{prop}
\vskip2mm
\pf~ Use Lemma \ref{lem2.9} first  to  infer  that for $0<|x|<N$,
 $$P_x[ \tau_{(-\infty, -N]} \vee  \tau_{[N,\infty)}<\tau_0]\leq C\frac{|x|}{N}\sup_{z>N} \Big(P_z[\tau_{(-\infty, -N]}<\tau_0]+P_{-z}[\tau_{[N,\infty)}<\tau_0]\Big)=o\bigg(\frac{x}{N}\bigg);$$
and  then, by employing  the inclusion-exclusion formula, to obtain 
the relation of the lemma. ~~ \qed

\vskip4mm\n
{\bf 2.3.}   ~In the following two lemmas we suppose that the walk $S_n$ is  aperiodic (i.e.,  $d_\circ =1$).

\begin{lem}\label{lem2.6}~ 
Let $d_\circ=1$. Then uniformly for $x, y\in \Z$, as $n\to\infty$
\begin{eqnarray}\label{eq2.6}
&&p^n(y-x)-p^n(-x)-p^n(y)+p^n(0) \nonumber\\
&&={\sf g}_{n}(y-x)-{\sf g}_{n}(-x)-{\sf g}_{n}(y)+{\sf g}_{n}(0)+o({xy}{n^{-3/2}}). ~~~~~~~~~~~~
\end{eqnarray}
\end{lem}
\v2\n
\pf~~Let  $\phi(l)$  denote  the characteristic  function of $Y$: 
$\phi(l)=E e^{ilY},$ $l\in \R$. As in  the usual proof of the local central limit theorem  choose a positive constant $\e$ so that $|\phi(l)-1|\geq \sigma^2 l^2/4$ for $|l|<\e$ and set $\eta=\sup_{\e\leq |l|\leq \pi}|1-\phi(l)|<1$.  Then the error in (\ref{eq2.6}) that we are to show to be $o(xyn^{-3/2})$  is written as
\beq
 \,(2\pi)^{-1/2}\int_{-\e}^\e \Big([\phi(l)]^n-e^{-n_* l^2/2}\Big)K_{x,y} (l)dl+O(e^{-n_*\e^2/2}+\eta^n)
\eeq
where 
$K_{x,y} (l)=e^{-i(y-x)l}-e^{ixl}-e^{-iyl}+1.$ Since $K_{x,y} (l)=(e^{ixl}-1)(e^{-iyl}-1)$, we have $|K_{x,y} (l)|\leq |xy|l^2$ and, scaling $l$ by $\sqrt{n_*}$ and applying the dominated convergence theorem,  we deduce that the integral above is $o(xyn^{-3/2})$ as required. ~~~\qed

\begin{lem}\label{lem2.7}~  Let $d_\circ=1$. Then uniformly in $y\in \Z$,    
$p^n(y)-p^n(0)={\sf g}_{n}(y)- {\sf g}_{n}(0)+o(y/{n})$ as $n\to\infty$; in particular $|p^n(y)-p^n(0)|\leq C|y|/{n}$.
\end{lem}
\v2\n
\pf~ The proof is similar to the preceding one. We have only to use  $|1-e^{ixl}|\leq |xl|$ in place of the bound of $K_{x,y}(l)$. ~~\qed

\section{Estimation of $q^n(x,y)$}

In this section we prove Theorem \ref{thm1.1}.   The proof relies on  the asymptotic estimate of the hitting-time distribution
$$f_x^{\{0\}}(k)= P_x[\tau_{0}=k]~~~~~(k=1,2,\ldots)$$
as $k\to\infty$, where $\tau_0$ denotes, as in Section 2, the first  time that $S^x_n$ hits the origin after  time 0. The following theorem is essentially proved in \cite{U}.
\v2\n
{\bf Theorem A} ~ ~{\it Under the basic assumption of this paper,   as $|x|\vee k\to \infty$}
\begin{eqnarray}\label{eqA}
f_x^{\{0\}}(k)&=&\frac{\sigma a^*(x)e^{- x^2/2\sigma^2k}}{\sqrt {2\pi}\, k^{3/2}}
+o\bigg( \frac{|x|+1}{k^{3/2}}\wedge \frac{1}{|x|^{2}+1}\bigg).
\end{eqnarray}

\v2\n
\pf  Immediate from Theorems 1.1 and 1.2 of \cite{U}.~~~ \qed
\v2
We have only to  consider  the case $0< |y|\leq x$ in view of the duality of $q^n(x,y)$ and $q^n(y,x)$ (i.e.,  transformed  to each other by  time-reversion). 

Let  $\phi(l)$  denote  the characteristic  function of $Y$ as in the proof of Lemma \ref{lem2.6}.
In what follows we suppose that the walk $S_n$ is aperiodic so that $|\phi(l)|<1$ for $|l|\leq \pi$.  We shall use the representation
\beqn\label{eq3.1}
q^n(x,y)=p^n(y-x)-\sum_{k=1}^n f_x^{\{0\}}(n-k)p^k(y)
\eeqn
and its Fourier version
\beqn\label{eq3.2}
q^n(x,y)=\frac1{2\pi}\int_{-\pi}^\pi \Big[\pi_{y-x}(t)-\rho(t)\pi_{-x}(t)\pi_y(t)\Big]e^{-int}dt~~~~~(x\neq 0)
\eeqn
and
\beqn\label{eq3.20}
q^n(0,y)=\frac1{2\pi}\int_{-\pi}^\pi \rho(t)\pi_y(t)e^{-int}dt,
\eeqn
where 
$$\pi_x(t)=\lim_{r\uparrow 1}\sum_{n=0}^\infty p^n(x)e^{itn}r^n=\frac1{2\pi}\int_{-\pi}^\pi \frac{e^{-ix l}}{1-e^{it}\phi(l)}dl~~~~~~(t\neq 0)$$
 and $\rho(t)=1/\pi_0(t)$; it holds that
 $$\rho(t)= \sigma\sqrt{-2it}(1+ o(1))~~~~~~\mbox{as} ~~~t\to 0$$
  (cf. \cite{U},  Section 2). Note that $q^n(0,y)=f^{\{0\}}_{-y}(n)$ by duality (or by coincidence of the Fourier coefficients), so that Theorem \ref{thm1.1} in the case $x=0$ is immediate from Theorem A.

  The supposition  that the  walk $S_n$ is  aperiodic gives rise to no essential loss of generality. To see  this 
let $d_\circ>1$ and put $\om=2\pi/d_\circ$. Then one can find a number $\xi$ among $1,\ldots, d_\circ-1$ such that $p(x+\xi)=0$ for all $x\notin d_\circ\Z$.   Hence  for all  $l\in (-\pi,\pi]$,
$$\phi(l+\om)=\sum e^{i(x+\xi)(l+\om)}p(x+\xi)=e^{i\xi\om}\phi(l).$$
Owing to the irreducibility of the walk there exists an  integer $k$ such that $k\xi =1$(mod$(d_\circ))$. Noting $\phi(\l+k\om)=e^{i\om}\phi(l)$,  one observes  that 
 $$\pi_x(t- \om)=\frac1{2\pi}\int_{-\pi}^\pi\frac{e^{-ix(l+k \om)}}{1-e^{it}e^{-i\om}\phi(l+k\om)}dl=\pi_x(t)e^{-ixk\om};$$
  in particular $\rho(t-\om)=\rho(t)$. It accordingly  follows that the integrand of the integral on the right side of (\ref{eq3.2}) is invariant by a shift of $t$ by $\om$ if (and only if) $(y-x)k\om -n\om\in 2\pi \Z$, namely $(y-x)k = n$(mod$(d_\circ))$ (the only case when $q^n(x,y)\neq 0$), hence the general case is reduced to the case $d_\circ=1$ since all our estimation of $q^n(x,y)$ is based on (\ref{eq3.2}).

\begin{thm}\label{thm3.1} Uniformly for $0<|y|\leq x< a_\circ \sqrt n$, as $n\to\infty$ and $|y|/\sqrt n\to 0$ 
$$q^n(x,y)={\sf g}_{n}(x)\frac{\sigma^4a(x)a(-y)+xy}{n_*}+o\bigg(\frac{xy}{n^{3/2}}\bigg).$$
\end{thm}
\v2
\pf First consider the case when not only $y$ but also $x$ is $o(\sqrt n)$. Of the integrand in
 (\ref{eq3.2}) make  the decomposition
\beq
\pi_{y-x}(t)-\rho(t)\pi_{-x}(t)\pi_y(t)&=&\pi_{y-x}-\pi_{-x}-\pi_{y}+\pi_{0}
+a(x)a(-y)\rho\\
&&-\, \rho\,{\rm e}_{x}\,{\rm e}_{-y}+a(x)\rho\,{\rm e}_{-y}+a(-y)\rho\,{\rm e}_{x},
\eeq
where
$$\,{\rm e}_{x}=\,{\rm e}_x(t)=\pi_{-x}(t)-\pi_0(t)+a(x).$$
Noting that $e^{-(\xi-\eta)^2}-e^{-\xi^2}-e^{-\eta^2}+1=e^{-\xi^2-\eta^2}(e^{2\xi\eta}-1)+O(\xi^2\eta^2)=2\xi\eta+o(\xi\eta)$ as $\xi, \eta\to 0$, we apply  Theorem A and Lemma \ref{lem2.6} to see
\beq
&&\frac1{2\pi}\int_{-\pi}^\pi \Big[\pi_{y-x}-\pi_{-x}-\pi_{y}+\pi_{0}
+a(x)a(-y)\rho\Big]e^{-int}dt\\
&&=p^n(y-x)-p^n(-x)-p^n(y)+p^n(0)+a(x)a(-y)f_0^{\{0\}}(n)\\
&&={\sf g}_{n}(0)\frac{\sigma^4a(x)a(-y)+xy}{n_*}+o\bigg(\frac{xy}{n^{3/2}}\bigg).
\eeq
In \cite{U} (Section 3) we have made decomposition $(2\pi)\,{\rm e}_x(t)=\,{\rm c}_x(t)+i\,{\rm s}_x(t)$, where
$${\rm c}_x(t)=\int_{-\pi}^\pi \bigg(\frac1{1-e^{it}\phi(l)}-\frac1{1-\phi(l)}\bigg)(\cos xl -1)dl$$
$${\rm s}_x(t)=\int_{-\pi}^\pi \bigg(\frac1{1-e^{it}\phi(l)}-\frac1{1-\phi(l)}\bigg)\sin xl\,dl$$ 
and verified the estimates given in the following two lemmas.
\v2\n
{\bf Lemma B1}~ {\it There exists a constant $C$ such that}
$$~~~~|{\rm c}_x(t)|\leq Cx^2\sqrt{|t|}, ~~~|{\rm c}'_x(t)|\leq  C x^2/\sqrt{|t|},~~~|{\rm c}''_x(t)|\leq Cx^2/|t|^{3/2}.$$

\vskip2mm\n
{\bf Lemma B2} ~~{\it Suppose that $E |Y|^{2+\de}<\infty$ for some $0\leq \de<1$.   Then, uniformly in $x\neq 0$, as $t\to 0$}
$$~~~~|{\rm s}_x(t)|/|x|=o(|t|^{\de/2}), ~~~|{\rm s}'_x(t)|/|x|=o(|t|^{\de/2}/|t|), ~~~|{\rm s}''_x(t)|/|x|=o(|t|^{\de/2}/|t|^2).$$

\v2
By a simple change of variables we derive the bounds
\beqn\label{114}
|\pi_{-x}^{(j)}(t)| \leq C|t|^{-\frac12 -j}~~~~~~~(j=0, 1, 2),
\eeqn
which in particular give  $|\rho^{(j)}(t)|\leq C|t|^{\frac12 -j}$,  where  the super script $(j)$ indicates the  derivative of $j$-th order.
With the help of  the  bounds given above as well as of Lemmas B1 and B2   we can readily infer   that each  of the contributions of $\rho\,{\rm e}_{x}\,{\rm e}_{-y}$, $a(x)\rho\,{\rm e}_{-y}$, and $a(-y)\rho\,{\rm e}_{x}$  to the integral in (\ref{eq3.2})  is $o(xy/n^{3/2})$. Eg., writing $g(t)=(\rho\,{\rm c}_{x}\,{\rm c}_{-y})(t)$, integrating by parts and observing that $|g'(t)|\leq Cx^2y^2 \sqrt{|t|}$ we obtain 
\beqn\label{lem}
\int_{-\pi}^\pi  g(t)e^{-int}dt= \frac{1}{in}\int_{-\pi}^\pi g'(t) e^{-int}dt= O\bigg(\frac{x^2y^2}{n^2\sqrt n}\bigg)+\int_{1/n<|t|<\pi} g'(t) e^{-int}dt.
\eeqn
Integrate by parts once more and apply the bound $|g''(t)|\leq Cx^2y^2/ \sqrt{|t|}$ to evaluate   the last integral to be $O(x^2y^2/n^2\sqrt n)$, which is  $o(xy/n^{3/2})$  since $x=o(\sqrt n\,)$.

It remains to consider the case $\e\sqrt n <x<a_\circ \sqrt n$. This time we use the decomposition
$$
\pi_{y-x}-\rho\pi_{-x}\pi_y=\pi_{y-x}-\pi_{-x}+a(-y)\rho\,{\pi}_{-x}+\, \rho\,{\pi}_{-x}\,{\rm e}_{-y}.
$$
Owing to Lemma \ref{lem2.7} and the present assumption on $x, y$, $p^n(y)-p^n(0)=o(y/n)=o(xy/n^{3/2})$.
Hence, again by Lemma \ref{lem2.6}  and Theorem A,
\beq
\frac1{2\pi}\int_{-\pi}^\pi \Big[\pi_{y-x}-\pi_{-x}
+a(-y)\rho\,{\pi}_{-x}\Big]e^{-int}dt
&=&p^n(y-x)-p^n(-x)+a(-y)f_x^{\{0\}}(n)\\
&=&{\sf g}_{n}(x)\frac{\sigma^4a(x)a(-y)+xy}{n_*}+o\bigg(\frac{xy}{n^{3/2}}\bigg).
\eeq
 In the same 
way  as is argued at (\ref{lem}) the contribution of $\rho\,{\pi}_{-x}\,{\rm e}_{-y}$ to the integral in (\ref{eq3.2}) can be evaluated to be  $O(y^2/n^{3/2})+o(y/n)$, which is $o(xy/n^{3/2})$.

 Theorem \ref{thm3.1} has been proved.~~~\qed

\v2
Theorem \ref{thm3.1} implies  {\bf (i)} of Theorem \ref{thm1.1} (in view of the local central limit theorem). 
\v2
\begin{prop}\label{prop3.1}~ Uniformly for $x, y$ such that both $x$ and $|y|$ are between $ a_\circ^{-1}\sqrt n$ and  $a_\circ \sqrt n$, as $n\to\infty$ 
\beq
q^n(x,y)&=&{\sf g}_{n}(y-x)-{\sf g}_{n}(y+x)+o({1}/{\sqrt n}) ~~~~~~\mbox{if}~~~~~y>0,\\
&=&o({1}/{\sqrt n})  ~~~~~~~~~~~~~~~~~~~~~~~~~~~~~~~~~~~\mbox{if}~~~~~y<0.
\eeq
\end{prop}
\v2\n
\pf  We  prove the second relation first.  To this end we introduce an auxiliary 
 walk. 
Let $\tilde p(x)$ be any probability law on $\Z$ of zero mean and  variance $\sigma^2$ such that its third absolute moment is finite and the random walk  determined by $\tilde p$ is left continuous, namely $\tilde p(y)=0$ for $y\leq - 2$. Let $\tilde p^n(y-x)$ and $\tilde q^n(x,y)$ denote the corresponding $n$-th step transition probabilities and let $y<0<x$.  From the assumed left continuity it follows that   $\tilde q^n(x,y)=0$, and  hence
\begin{eqnarray}\label{I-III}
q^n(x,y)&=& p^n(y-x)-\tilde p^n(y-x)-\sum_{k=1}^n \Big(f_x^{\{0\}}(k)p^{n-k}(y)-\tilde f_x^{\{0\}}(k)\tilde p^{n-k}(y)\Big) \nonumber \\
&=& U_n(x,y)+V_n(x,y)+W_n(x,y),
\end{eqnarray}
where
$$U_n(x,y)=p^n(y-x)-\tilde p^n(y-x), $$
$$
V_n(x,y)=-\sum_{k=1}^{n-1} \Big( f_x^{\{0\}}(k)- \tilde f_x^{\{0\}}(k)\Big) p^{n-k}(y)
$$
and
$$W_n(x,y)=-\sum_{k=1}^{n-1} \tilde f_x^{\{0\}}(k)\Big( p^{n-k}(y)- \tilde p^{n-k}(y)\Big).$$

\v2\n
By the local limit theorem $U_n=o(1/\sqrt n)$.  We apply Theorem A to see that $\sup_{k\ge 1}|f_x^{\{0\}}(k)-\tilde f_x^{\{0\}}(k)|=o(1/x^2)=o(1/n)$, which combined with the  trite bound   $\sup_z p^k(z)\leq C/\sqrt k$ shows $V_n=o(1/\sqrt n)$. The bound  $W_n=o(1/\sqrt n)$ is verified e.g. by observing that $\sup_{1\leq k\le n}|p^k(y)-\tilde p^k(y)|\sqrt k\to 0$ as $n\to \infty$ uniformly in $y$. 
For the proof of the first relation  we write for $y>0$
\beq
q^n(x,y)&=&p^n(y-x)-\sum_{k=1}^n f_{-x}^{\{0\}}(k)p^{n-k}(y)+ \sum_{k=1}^n [f_{-x}^{\{0\}}(k)- f_x^{\{0\}}(k)]p^{n-k}(y)\\
&=&p^n(y-x)-p^n(y+x) +q^n(-x,y)+r_n(x,y),
\eeq
where $r_n(x,y)=\sum_{k=1}^n [f_{-x}^{\{0\}}(k)- f_x^{\{0\}}(k)]p^{n-k}(y)$.
In view of what has been shown above as well as the local central limit theorem it suffices to show $r_n(x,y)=o(1/\sqrt n)$. There exists a positive integer $N$ such that for $0<\e<1/2$ and $n>N$,
$$\sum_{1\le k<\e n} f_{\pm x}^{\{0\}}(k)p^{n-k}(y)\le  P_{\pm x}[\tau_0 <\e n]/\sqrt{\pi n_*}$$
and in view of Donsker's invariance principle the probability on the right side above tends to zero as $\e\downarrow 0$ uniformly for $x>\sqrt n/a_\circ$.  Now the required estimate follows from Theorem A, according to which 
  $f_{-x}^{\{0\}}(k)- f_x^{\{0\}}(k)=o(x/k^{3/2})$ as $x\wedge k\to\infty$. ~~\qed

\v2\n

\v2
\begin{prop}\label{prop3.2}~ Suppose $E|Y|^{2+\de}<\infty$ for some $ \de\geq 0$. Then,  uniformly for
 $ |x| <a_\circ\sqrt n$ and $ |y|> a_\circ^{-1} \sqrt n$, as $n\to\infty$
$$q^n(x,y)=O\bigg(\frac{x}{y}{\sf g}_{4n}(y)\bigg)+o\bigg(\frac{x}{|y|^{2+\de}}\bigg).
 $$
\end{prop}
\v2\n
\pf     Suppose $y/2>\sqrt {n_*}$ for simplicity. ~Put $\tau=\tau_0\wedge \tau_{(y/4,\infty)}$. Then
\beq
q^n(x,y)&=&P_x[\tau\le n <\tau_0, S_n=y]\\
&=&P_x[y/4<S_\tau<y/2,  n <\tau_0, S_n=y]+P_x[S_{\tau} \ge y/2, n <\tau_0, S_n=y]\\
&=& I+ II~~~~\mbox{(say).}
\eeq 
We employ the inequality
$$I\leq \sum_{k=1}^n\sum_{y/4<z<y/2} P_x[\tau=n-k, S_\tau=z]P_z[S_k=y].$$
The following less familiar version of local central limit theorem is found in \cite{U0} (see its Corollary 6): under the assumption of Proposition \ref{prop3.2}
\beqn\label{eqLLT}
P_0[S_n=x]
={\sf g}_n(x)\left[1+ P^{n,\nu}(x) \right] 
+o\left(\frac1{\sqrt n^{1+\de}}\wedge\frac{\sqrt n }{|x|^{2+\de}}\right),
\eeqn
 $(n+|x|\to \infty)$, where $\nu=\lfloor \de \rfloor$ (the largest integer that  does not exceeds $\de$), $P^{n,0}\equiv 0$ and 
$P^{n,\nu}(x)= \frac1{\sqrt n} P_1\left(\frac{x}{\sqrt n}\right)+\cdots +\frac1{\sqrt n^{\,\nu}} P_{\nu}\left(\frac{x}{\sqrt n}\right)$ if $\nu\geq 1$ with  the same real polynomials $P_j$ of degree $j$ as those  associated with the Edgeworth expansion.
From (\ref{eqLLT}) one deduces
$$\max_{1\leq k\leq n}\, \max_{y/4<z<y/2}P_z[S_k=y]=O\bigg({\sf g}_{4n}(y)\bigg)+o\bigg(\frac{\sqrt n}{y^{2+\de}}\bigg)$$
(use $(y/2)^2> n_*$ for evaluation of the maximum over $k$). On the other hand
$$\sum_{k=1}^n\sum_{y/4<z<y/2} P_x[\tau=n-k, S_\tau=z]\le P_x[\tau_{(y/4,\infty)}<\tau_0]=O\bigg(\frac{x}{y} \bigg).$$
Hence
\beqn\label{eqI}
I= O\bigg(\frac{x}{y}{\sf g}_{4n}(y)\bigg)+o\bigg(\frac{x\sqrt n}{|y|^{3+\de}}\bigg).
\eeqn

For  evaluation of  $II$ we begin with
\[
II\leq \sum_{k=1}^n P_x[S_{k} \ge y/2, \tau=k, S_n=y].
\]
Under $S_k\geq y/2$ we have $\{\tau=k\}=\{\tau>k-1\}$, hence the sum on the right side equals
$$ \sum_{k=1}^n E_x\Big[P_{S_{k-1}}[S_{1} \ge y/2, S_{n-k+1}=y];~ \tau >k-1\, \Big].
$$ 
Since  $P_x[\tau> k-1]\le P_x[\tau_0> k-1]=O(x/\sqrt k\,)$ and for $z<y/4$,
$$P_{z}[S_{1} \ge y/2, S_{n-k+1}=y]\leq \sum_{w\geq y/2}p(w-z)p^{n-k}(y-w)= o\bigg(\frac{1}{y^{2+\de}\sqrt{n-k+1}}\bigg),$$
we get
$$II=\sum_{k=1}^n P_x[\tau>k-1]\times o\bigg(\frac{1}{y^{2+\de}\sqrt{n-k+1}}\bigg)=
o\bigg(\frac{x}{y^{2+\de}}\bigg).
$$
This together with (\ref{eqI}) shows the  estimate of the proposition. ~~~\qed


\section{Estimation of $q^n_{(-\infty,0]}$ and $h_x(n,y)$}

{\it Proof of Proposition \ref{thm1.2}.}~ In view of {\bf (i)} and {\bf (ii)} of Theorem \ref{thm1.1}  it suffices to 
prove that uniformly  in $n$, 
$$q^n(x,y)- q_{(-\infty,0]}^n(x,y)= o\bigg(\frac{xy}{n^{3/2}}\bigg)~~~~~~~\mbox{as}~~~~~~~x\wedge y\to\infty.$$
This difference  may be written as
\beqn\label{4.1}
\sum_{k=1}^n\sum_{z<0} h_x(k,z) q^{n-k}(z,y).
\eeqn
  Employing  the identity (\ref{eq2.51}) one observes that
$$\sum_{1\leq k<n/2}\,\sum_{z<0} h_x(k,z)|z|\leq  \sum_{z<0} H_x^+(z)|z|=f_+(x)-x=o(x)~~~~\mbox{as}~~~x\to\infty.$$
Combined with the simple bound (\ref{iv}) this shows that  the sum over $k\leq n/2$ in (\ref{4.1}) is $o(xy/n^{3/2})$. 
 The other half of the sum  is less than the probability that the time-reversed walk starting at $y$ enters $(-\infty,0]$ till the time $n/2$ and ends in $x$ at the time $n$ and hence estimated  also to be $o(xy/n^{3/2})$. ~~~\qed
 \v2\n
 \begin{lem}\label{lem4.1}~~For each $x=1, 2, \ldots$, uniformly for $n\geq y^2/a_\circ$, as $y\to\infty$
$$q_{(-\infty,0]}^n(x,y)=\frac{2f_+(x)y}{n_*}{\sf g}_n(y)(1+o(1)).$$
\end{lem}
\vskip2mm\n
{\it Proof.}~ Given a positive integer $x$, take an integer $N>x$ and put $\tau=T\wedge \tau_{[N,\infty)}$, the first leaving time from $[1, N-1]$. Then 
\beqn\label{eq4.0}
q_{(-\infty,0]}^n(x,y)= E_{x}[ q_{(-\infty,0]}^{n-\tau}(S_{\tau},y); \tau<T\wedge (n+1)] .
\eeqn
Let $\a$ be any positive number less than 1.   For each $\e>0$ we can choose  $N$ large enough  that  for all $ k, n, z$ and $y$ that satisfy $0\leq k< n^{\a}$,  $2N<y\leq \sqrt{a_\circ  n}$ and $ N\leq z\leq \sqrt{n}/N$,  the following three
 bounds hold:
\beqn\label{eq4.1}
\bigg|q_{(-\infty,0]}^{n-k}(z,y)-\frac{2zy}{n_*}{\sf g}_n(y)\bigg|< \frac{\e zy}{n_*^{3/2}},
\eeqn
\beqn\label{eq4.2}
|P_{x}[\tau<T\,]-f_+(x)/N|\leq \e x/N,
\eeqn
\beqn\label{eq4.21}
E_{x}[ S_{\tau}-N; \tau<T] \leq \e x,
\eeqn
according  to {\bf (i)} of Theorem \ref{thm1.1} and Proposition \ref{thm1.2} for (\ref{eq4.1}), to Proposition \ref{lem2.3} for (\ref{eq4.2}) and  to Proposition \ref{lem2.4} for (\ref{eq4.21}).
Since $\tau$ equals the sum of the sojourn times of  sites $w$ in the interval  $[1,N-1]$  spent by the walk before leaving it, we have $E_{x}[\tau] =\sum_{w=1}^{N-1}\sum_{k=0}^\infty P_x[ S_k=w, k<\tau]\leq \sum_{w=1}^{N-1}g_{(-\infty,0]}(x,w)\leq CxN$, and on using this
\begin{eqnarray}\label{ppp}
P_{x}[ S_{\tau}> \sqrt{n}/N, \tau< T]&\leq &\sum_{k=1}^\infty P_{x}[ S_{k}>\sqrt{n}/N, \tau = k] \nonumber \\
&\leq&\sum_{k=1}^\infty P_{x}[ Y_{k}> \sqrt{n}/N -S_{k-1},\, \tau > k-1]  \nonumber \\
&\leq& \sum_{k=0}^{\infty}P_{x}[\tau > k]P[Y>\sqrt{n}/N -N] = xN^3\times o\bigg( \frac{1}{n}\bigg)  ~~~~~~
\end{eqnarray}
as $n \to\infty$. Since $ q_{(-\infty,0]}^{n-k}(\cdot, y)\leq C'/\sqrt n$ if $k<n^\a$, this entails that
\beqn\label{eq4.3}
E_{x}[ q_{(-\infty,0]}^{n-\tau}(S_{\tau},y);S_{\tau}>\sqrt{n}/N, \tau<T\wedge n^\a\,] =o(n^{-3/2})
\eeqn
as $n\to \infty$ (with $N, x$ fixed).
On the other hand, using (\ref{eq4.1}) we obtain   
\begin{eqnarray}\label{pppp}
&&E_{x}[q_{(-\infty,0]}^{n-\tau}(S_{\tau},y);S_{\tau}\leq \sqrt{n}/N, \tau<T\wedge n^\a\, ]   \nonumber\\
&&=E_{x}\bigg[\frac{2Ny}{n_*}{\sf g}_n(y); S_{\tau}\leq \frac{\sqrt{n}}{N}, \tau<T\wedge n^\a\,\bigg]
 + E_{x}\bigg[\frac{2(S_{\tau}-N)y}{n_*}{\sf g}_n(y);S_{\tau}\leq \frac{\sqrt{n}}{N}, \tau<T\wedge n^\a\, \bigg]   \nonumber \\
&&~~+r(n,x,y)
\end{eqnarray}
 with 
$|r(n,x,y)|\leq \e y E_x[S_\tau ; \tau<T]/n_*^{3/2}$.  Writing $E_x[S_\tau ; \tau<T] = NP_x[\tau<T] +E_x[S_\tau -N; \tau<T]$  we apply  (\ref{eq4.21}) and (\ref{eq4.2}) to see that the remainder $r$ as well as the second expectation on the right side in (\ref{pppp}) is  dominated in absolute value by  $3\e xy/n_*^{3/2}$.
We also have $P_x[\tau> n^\a] = O(e^{-\k n^\a})$ with some $\k=\k_N>0,$ and hence,  owing to (\ref{ppp}),
$$P_x[S_\tau\leq \sqrt{n}/N, \tau< T\wedge n^\a]=P_x[\tau<T]+o(1/n).$$
 Combining these bounds with  
 (\ref{eq4.3}) and (\ref{eq4.0})  shows that for all sufficiently large $y$ and for $n>y^2/a_\circ$,
$$\bigg|q_{(-\infty,0]}^n(x,y)-\frac{2Ny}{n_*}{\sf g}_n(y)P_{x}[\tau<T\,]\bigg|<
\frac{7\e xy}{n_*^{3/2}}$$
and substitution from   (\ref{eq4.2})  completes the proof of Lemma \ref{lem4.1}.~~~\qed 

\begin{lem}\label{lem4.2}~~For each $x, y=1, 2, \ldots$, as $n\to\infty$
$$q_{(-\infty,0]}^n(x,y)=\frac{2f_+(x)f_-(y)}{n_*}{\sf g}_n(0)(1+o(1)).$$
\end{lem}
\vskip2mm\n
{\it Proof.}~ 
Applying Lemma \ref{lem4.1} to the time-reversed walk we have
$$
\bigg|q_{(-\infty,0]}^n(z,y)-\frac{2z{\sf g}_n(z)f_-(y)}{n_*}\bigg|< \frac{\e z{y}}{n^{3/2}}
$$
(valid for all  $z\geq N$) in place of (\ref{eq4.1}) and we can proceed as in the proof of Lemma \ref{lem4.1}.~~~\qed
 \v2
Theorem \ref{thm1.3}  follows from Proposition \ref{thm1.2} and Lemmas \ref{lem4.1} and \ref{lem4.2} given above.
\v2\v2\n
{\it Proof of Theorem \ref{thm1.4}.} ~ The probability $h_x(n,y)$ is represented as
\beqn\label{upb-h}
h_x(n,y)=  \sum_{z>0} q_{(-\infty,0]}^{n-1}(x,z)p(y-z).
\eeqn
Write $F(x,n)=2f_+(x){\sf g}_n(x)/n_*$.    In view of Theorem \ref{thm1.1} and a local limit theorem, for each $\e>0$ we can then choose $\eta>0$ such that for all sufficiently large $n$,
$$\Big|q_{(-\infty,0]}^{n-1}(x,z)-F(x,n)f_-(z)\Big|\le \e F(x,n)f_-(z)$$
whenever $0<z\le \eta \sqrt n$ and $0< x<a_\circ\sqrt n$.    Hence,  on using the second expression of $H_{\infty}^+$ in (\ref{q}), the difference  $|h_x(n,y)-F(x,n)H_{\infty}^+(y)|$ is at most
$$ \e F(x,n)H_{\infty}^+(y) +\sum_{z>\eta \sqrt n} \Big|q_{(-\infty,0]}^{n-1}(x,z)-F(x,n)f_-(z)\Big|p(y-z).$$
Owing to (\ref{iv})  the summand of the last sum is at most a constant multiple of $n^{-3/2}xzp(y-z)$, so that if  $\a_n(y)=\sum_{z> \eta \sqrt n}zp(y-z)$, then this sum is at most $n^{-3/2}x\a_n(y)$, hence
$$|h_x(n,y)-F(x,n)H_{\infty}^+(y)|\le \e F(x,n)H_{\infty}^+(y)+n^{-3/2}x\a_n(y).$$
 The proof is now finished by observing that if $E[|Y|^{2+\de};Y<0]<\infty$, then
$$\a_n(y)=o\bigg(\frac1{(\sqrt n +|y|)^{1+\de}}\bigg),~~~  \sum_{y\leq 0}\a_n(y)=o(n^{-\de/2})~~~\mbox{  and}~~~ \sum_{y\leq 0}\a_n(y)|y|^\de=o(1).$$
   The first half of Theorem \ref{thm1.4} has been verified.

For  the second half we verify that for $y\leq 0$ and $x\geq \sqrt n$,
\beqn\label{upb-h2}
\sum_{z=1}^\infty q^{n}(x,z)p(y-z)\leq C\bigg[\frac{{\sf g}_{4n}(x)}{n^{1/2}}+o\bigg(\frac1{x^{2+\de}}\bigg)\bigg]H_{\infty}^+(y) +\frac{C}{n^{1/2}}P[Y<y-{\textstyle \frac12}x].
\eeqn
Since $h_x(n,y)$ is not larger than the sum on  the left side, this implies (\ref{upb-h20}).  For verification of (\ref{upb-h2}) we
break the range of summation  into three parts $0<z\leq \sqrt n \wedge \frac12 x$,  $\sqrt n <z\leq x/2$ and $z>x/2$,
and denote the corresponding sums by $I$, $II$ and $I\!I\!I$, respectively.
It is immediate  from {\bf (iii)} of Theorem \ref{thm1.1} that  $I =(O( {\sf g}_{4n}(x)x^{-1}+o(x^{-2-\de}))H_{\infty}^+(y)$.  The local limit theorem estimate (\ref{eqLLT})  gives that $q^{n}(x,z)\leq p^n(z-x)\leq 2{\sf g}_{n}(x/2)+ o(n^{1/2}x^{-2-\de})$ for  $\sqrt n<z\leq x/2$ and   we apply this bound as well as the  bound
\beqn\label{eq4.-1} \sum_{z>\sqrt n}p(y-z) \leq \frac1{\sqrt n}\sum_{z=1}^\infty zp(y-z) \leq \frac2{\sigma^2\sqrt n}\bigg[\sup_{z\geq 1}\frac{f_-(z)}{z}\bigg] H_{\infty}^+(y)
\eeqn
 to have $II=[O( {\sf g}_{4n}(x)n^{-1/2}) +o(x^{-2-\de})]H_{\infty}^+(y)$. Finally 
$I\!I\!I \leq C n^{-1/2}\sum_{z<y-x/2}p(z)$.
 These
estimates together verify  (\ref{upb-h2}). As in (\ref{eq4.-1})  we derive $P[Y<y-{\textstyle \frac12}x]\leq C_1 H^+_{\infty}(y)/x$. We also have ${\sf g}_{4n}(x)\leq C_1/x$. Hence the last relation of the theorem follows from  (\ref{upb-h20}). The proof of Theorem \ref{thm1.4} is complete.
~~~\qed

\v2\v2
  The proof  of Theorem \ref{thm1.5} is based on Theorem \ref{thm1.4}.    Put
  $$\Phi_\xi(t)=\frac{|\xi| e^{-\xi^2/2t}}{\sqrt{2\pi}\, t^{3/2}}~~~~~~~(t>0, \xi \neq 0).$$
 Then  Theorem \ref{thm1.5} (under $d_\circ =1$)  may be  restated as follows.
\v2\n
{\bf Theorem 1.2} ~~{\it Suppose that $E[|Y|^3;Y<0]<\infty$. Let $y<0<x$. Then uniformly for $x,|y|\le a_\circ\sqrt n$, as $n\to\infty$ and $x\wedge|y|\to \infty$}
\beqn\label{prop-est}
q^n(x,y)=C^+\Phi_{x+|y|}(n_*)+ o\bigg(\frac{x +|y|}{n^{3/2}}\bigg).
\eeqn

\v2\n
\pf ~ We use the representation
$$q^n(x,y)=\sum_{k=1}^n \sum_{z<0}h_x(k,z)q^{n-k}(z,y).$$
Break  the right side into three parts by partitioning the range of the first summation as follows
\beqn\label{eq4.5}
1\leq  k< \e n; ~~ \e n \leq k \leq (1-\e)n;~~ (1-\e)n<k\leq n
\eeqn
and call the corresponding sums $I,~II $ and $I\!I\!I$, respectively. Here $\e$ is a positive constant that will be chosen small. 

Consider the limit procedure as indicated in the theorem.  First suppose that $(x\wedge |y|)/\sqrt n$ is bounded away from zero. Then by (\ref{eq1.4}), the last relation of Theorem \ref{thm1.4}, and (\ref{iv})
$$I\leq  \sum_{1\leq k <\e n}\frac C {k^{1/2}\,x}\sum_{z<0}  H_{\infty}^+(z)\frac {|zy|}{n^{3/2}} \leq C' \frac{\sqrt{\e}}n.$$
and similarly $I\!I\!I \leq C'' \sqrt{\e}/n$ (see also (\ref{upb-h})). By the first half of Theorem \ref{thm1.4} and (\ref{iv})
$$II=\sum_{ \e n\leq k \le (1-\e)n}\frac{f_+(x){\sf g}_{k}(x)}{k}\,\sum_{z=-x}^{-1}H_{\infty}^+(z) q^{n-k}(z,y)(1+o_\e(1))+o_\e\bigg( \frac{1}{n}\bigg).$$
Here (and  in the rest of the proof) the estimate indicated by $o_\e$ may depend on $\e$ but is  uniform in the  limit under  consideration once $\e$ is fixed.  We substitute  from  {\bf (i)} of Theorem \ref{thm1.1} for $q^{n-k}$ and observe,  on  replacing $f_+(x)$  and $a(-y)$, respectively,  by $x$ and  $-y/\sigma^2$, 
\[
II=\sum_{ \e n\leq k \le (1-\e)n}\frac{x|y|{\sf g}_{k}(x){\sf g}_{n-k}(y)}{\sigma^2 k(n-k)}\,\sum_{z=-x}^{-1}H_{\infty}^+(z)(\sigma^2a(z)-z)(1+o_\e(1))+o_\e\bigg( \frac{1}{n}\bigg) \nonumber.
\]
Noting $x{\sf g}_k(x)/ k= \Phi_{x/\sigma}(k)$, we see
$$\sum_{ \e n\leq k \le (1-\e)n}\frac{x|y|{\sf g}_{k}(x){\sf g}_{n-k}(y)}{\sigma^2 k(n-k)}=\frac1{n\sigma^2}\int_0^1 \Phi_{x/\sqrt{n_*}}(t)\Phi_{y/\sqrt{n_*}}(1-t)dt + O\bigg(\frac{\e}{n}\bigg) +o\bigg( \frac{1}{n}\bigg).$$
Here we have used the assumption that $(x\wedge |y|)/\sqrt n$ is bounded away from zero as well as from infinity.
Since  $\Phi_\xi$ is the density of a Brownian passage-time distribution, we have
$$\int_0^1 \Phi_{\xi}(t)\Phi_{\eta}(1-t)dt=\Phi_{|\xi|+|\eta|}(1).$$
Hence
\beqn\label{II}
II= \frac1{n\sigma^2} \Phi_{(x+|y|)/\sqrt{n_*}\,}(1)\sum_{z=-x}^{-1}H_{\infty}^+(z)(\sigma^2a(z)-z)+ O\bigg(\frac{\e}{n}\bigg)+ o_\e\bigg( \frac{1}{n}\bigg).
\eeqn
Recalling that $C^+=\sum _{z<0}H_{\infty}^+(z)(\sigma^2 a(z)-z)$ we then see $\sigma^2 nII-C^+ \Phi_{(x+|y|)/\sqrt{n_*}}(1) \to 0$  (as well as  $nI+n I\!I\!I \to 0$)  as $n\to\infty$ and $\e\to 0$ in this order.  Thus  (\ref{prop-est}) is obtained.

Next suppose $x\wedge |y|= o(\sqrt n\,)$.
By duality one may suppose that $x=o(\sqrt n)$. 
From Theorem \ref{thm1.4} (with $\de=1$) and  from the bound $H^+_x(y)\leq CH_{\infty}^+(y)$ (see (\ref{17})) one deduces, respectively,  
\beqn\label{hh}
{\rm i)}~~~\sum_{k\geq \e n} \sum_{z<0} h_x(k,z)|z|\le M_\e x/\sqrt{ n}~~~~~\mbox{and}~~~~~{\rm ii)}~~~ \sum_{z<-x} H_x^+(z)z =o(1).~~~
\eeqn
Here (and below) $M_\e$ indicates  a constant that may depend on $\e$ but not on the other variables.   On using i) above with the help of the bound $q^k(z,y)\leq C|zy|k^{-3/2}$
$$II \leq M_\e xy/n^2= o_\e( {y}{n^{-3/2}})$$
(as $n\to\infty$  under the supposed constraints on $x, y$); 
similarly on using Theorem \ref{thm1.1} {\bf(i)} together with ii) above
$$I=\sum_{1\le k< \e n}\,\sum_{z=-x}^{-1}h_x(k,z)\cdot \frac{\sigma^4a(z)a(-y)+zy}{\sigma^2(n-k)}{\sf g}_{n-k}(y)(1+o_e(1))+o\bigg( \frac{y}{n^{3/2}}\bigg).$$
For the evaluation of the last double sum we may replace $(n-k)^{-1}{\sf g}_{n-k}$ by $n^{-1}{\sf g}_{n} (1+O(\e))$.  Since  $x\wedge|y|$ is supposed to go to infinity, we may also replace $a(-y)$  by  $|y|/\sigma^2$ and  in view of (\ref{hh}) we may extend the range of the double summation in the above expression of $I$   to  the whole quadrant $k\geq 1, z<0$;  moreover the sum  $\sum_{z=-\infty}^{-1}H^+_x(z)[\sigma^2 a(z)-z]$ that accordingly comes out and equals $\sigma^2 a(x)-x$ may be replaced by $C^+$ (see Lemma \ref{lem2.50}, Corollary \ref{lem2.5} and (\ref{eq2.51})). 
This leads to 
 $$ I=  C^+|y|{\sf g}_{n}(y)n_*^{-1}(1+O(\e))+o_\e(y{n^{-3/2}}).$$
As to $ I\!I\!I$   first observe that
 $$\sum_{k=1}^{\e n} q^k(z,y)=g_{\{0\}}(z, y)-r_n\leq C(|z|\wedge |y|)~~~\mbox{with}~~~0\leq r_n\leq C|zy|/\sqrt{\e n},$$
as follows from  (\ref{iv}) and (\ref{g}). 
 If $y/\sqrt n$ is bounded away from zero, then  $ I\!I\!I=O(x/n^{3/2})=o(y/n^{3/2})$. On the other hand, applying Theorem \ref{thm1.4} we  find that if $y=o(\sqrt n)$,
 $$ I\!I\!I= f_+(x){\sf g}_{n}(x)n_*^{-1}\sum_{z<0}H_{\infty}^+(z)g_{\{0\}}(z, y)(1+O(\e))+ o_\e\Big(x{n^{-3/2}}\Big),$$
hence in view of  $g_{\{0\}}(z, y)=a(z)-\sigma^{-2}z(1+o(1))$ (as $y\to -\infty$ uniformly for $z<0$)
  $$ I\!I\!I=  C^+x{\sf g}_{n}(x)n_*^{-1}(1+O(\e))+o_\e\Big(x{n^{-3/2}}\Big).$$
Adding these contributions yields the desired formula. ~~~\qed

\section{Estimation of $Q^+_n$}

{\it Proof of Proposition \ref{prop1.3.1}.}~ We apply  {\bf (i)} and {\bf (iii)} of Theorem \ref{thm1.1}.  On noting that $\sum_{y<0}(\sigma^2 a(-y)+y){\sf g}_n(y)=o(\sqrt n)$,  as $|x|/\sqrt n\to 0$, 
\[
Q_x^+(n) =  \frac{\sigma^2a(x)-x}{n_*}\sum_{-a_\circ\sqrt n<y<0} (-y){\sf g}_n(y)\Big[1+o(1)\Big]
+O\bigg(x\sum_{y\le -a_\circ\sqrt n}\frac {{\sf g}_{4n}(y)}{-y} \bigg) 
+o\bigg(\frac{x}{n^{1/2}}\bigg),
\]
which shows the first assertion of Proposition \ref{prop1.3.1} since $a_\circ$ is arbitrary and $\int_0^\infty ue^{-u^2/2}du=1$.  

For the second one, the case $x=o(\sqrt n)$ follows from what has just been proved.  In view of Proposition \ref{thm1.2} and  {\bf (i)} of Theorem \ref{thm1.1},  it  therefore suffices to show that 
$\sum_y q^n_{(-\infty, 0]}(x,y)=\int_{-x}^x {\sf g}_n(t)dt (1+o(1))$ uniformly for $|x|\geq  a_\circ^{-1}\sqrt n $,
which however follows from  Donsker's invariance principle together with the reflection principle.

The last assertion  of the proposition follows from Theorem \ref{thm1.5}, {\bf (iii)} of Theorem \ref{thm1.1} and the  
fact that $\int_0^{a_\circ} \Phi_{|\xi|+\eta}(1)d\eta\to (2\pi)^{-1/2}$ as $\xi\to 0, a_\circ \to\infty$. ~~~\qed
\v2

\begin{lem} \label{lem5.1}~ ~~Suppose  $E[|Y|^{3}; Y<0]<\infty$. Then uniformly for $x> a_\circ^{-1}\sqrt n$,  
\beqn\label{Q-4}
Q_x^+(n)\leq C{\sf g}_{4n}(x)+o\Big(\sqrt n /x^{3}\Big) +C\sum_{y<0}|y|P[Y<y-{\textstyle \frac12}x];
\eeqn
in particular $\sum_{x> M \sqrt n} Q^+_x(n)\to 0$ as $M\to\infty$ uniformly in $n$.
\end{lem}
\v2
\pf~   As in the proof of Theorem \ref{thm1.5}
we use the representation
\beqn\label{Q-2}
Q_x^+(n)=\sum_{k=1}^n \sum_{z<0}h_x(k,z)Q^+_z({n-k}),
\eeqn
By Proposition \ref{prop1.3.1} $Q_z^+(n-k)\leq C|z|/\sqrt{n-k}$ for all $z$ since $Q^+_x(n)\leq 1$.  From  the assumption  $E[|Y|^{3}; Y<0]<\infty$ it follows that $\sum_{z<0}H_\infty^+(z)|z|<\infty$. Note that $\sum_{k=1}^{n-1} k^{-1/2}(n-k)^{-1/2}$ is bounded.  Substitution from (\ref{upb-h20}) then leads to the first
 estimate (\ref{Q-4}). The second  relation  is immediate from it if one notes that $\sum_{x=0}^\infty\sum_{y<0}|y|P[Y<y-{\textstyle \frac12}x]$ is dominated by a constant multiple of  $E[|Y|^3; Y<0]$.
~~~\qed

\v2\n
{\it Proof of Theorem \ref{thm1.3.2}.}  ~Suppose $E[|Y|^{3}; Y<0]<\infty$.  Then in view of Lemma \ref{lem5.1} 
we have only to evaluate $\sum_{1\leq x\leq  M \sqrt n} Q^+_x(n)$ for each $M$. Now apply Theorem \ref{thm1.5} 
with the observation that   $\int_{\xi>0}d\xi \int_{\eta>0} \Phi_{\xi+\eta}(1) d\eta = (2\pi)^{-1/2}\int_0^\infty e^{-\xi^2/2}d\xi=1/2$, 
and you immediately find the first formula of the theorem.  If $E[|Y|^{3}; Y<0]=\infty$, one  has only to look at the relation (\ref{II}) which is valid without the third moment condition; in fact,   $\sum_{x=1}^\infty Q_x^+(n)$ is bounded below by the sum of $II$ in (\ref{II})  over $1\leq x<\sqrt n$, $-\sqrt n<y<0$ and the latter diverges to $+\infty$. ~~~\qed

\v2\n
{\it Proof of Corollary \ref{cor1.2} }.~  Lemma \ref{lem5.1} shows that if one writes
$$E[N_n(\ell)]=\sum_{1\leq  x\leq M\sqrt {n_*}}\,\,\sum_{-\ell\sqrt {n_*}\leq y\leq -1}m_n(x) q^n(x,y)+ \e_M(n),$$
then $ \e_M(n)\to 0$ as $M\to\infty$ uniformly in $n$.
According to Theorem \ref{thm1.5} and the assumption on $m_n(x)$ the double sum on the right side is asymptotically equal to
$$C^+\int_0^M d\xi\int_0^{\ell} \Phi_{\xi+\eta}(1)d\eta=C^+\frac1{\sqrt{2\pi}}\int_0^M\big(e^{-\xi^2/2}-e^{-(\ell+\xi)^2/2}\big)d \xi.$$
Since $M$ is arbitrary, we may let $M\to\infty$ to find the desired formula.~~\qed

\section{Absorption at the origin with probability $\a\in (0,1)$}
Let $\a\in (0,1)$  and consider the walk that is absorbed with probability $\a$ and continues to walk
with probability $1-\a$  every time when it is about to visit  the origin (thus the walk visits the origin if it is not absorbed, while it  does not and disappears  if it is). Let $q^n_\a(x,y)$ be the $n$-th step transition probability of it (set $q_\a^0(x,y)={\bf 1}(x=y)$ as usual) and denote by $r^n_{\a}(x,y)$ the probability that this walk starting at $x$ has visited the origin  but not been absorbed by the time $n$ when it is at $y$, so that
$$q_\a^n(x,y)=q^n(x,y)+ r_\a^n(x,y).$$
\begin{prop}\label{prop6.1}~~ Let $d_\circ=1$. Then uniformly for $|x|\vee |y|<a_\circ \sqrt n$,~  as $n\to\infty$ 
\beqn\label{rr}
r_\a^n(x,y)=\frac{1-\a}{\a}\cdot \frac{\sigma^2 [a^*(x)+ a^*(-y)]}{ n}\,{\sf g}_n(|x|+|y|)(1+o(1)).
\eeqn
\end{prop}
\v2\n
\pf~~
 Set  $f_x^{(1)}(k)=f_x^{\{0\}}(k)$ and, for $j=2, 3,\ldots$, inductively define 
$f_x^{(j)}(k) =\sum_{l=1}^k f_x^{(j-1)}(k-l)f_0^{(1)}(l)$
(the probability that the $j$-th visit of the origin occurs at time $k$). Then for $n=1,2,\ldots$,
$$r_\a^n(x,y)=\sum_{j=1}^\infty\sum_{k=1}^n(1-\a)^j f_x^{(j)}(k)q^{n-k}(0,y).$$ 
(This is valid even for  $y=0$ when the second sum concentrate on $k=n$.) 
We have  $\hat f_x^{(j)}(t):=\sum_{k} f_x^{(j)}(k) e^{ikt}=\hat f_x^{\{0\}}(t)[\hat f_0^{\{0\}}(t)]^{j-1}.$ One can readily derive
 \beqn\label{111}
 \hat f_x^{\{0\}}(t)=\pi_{-x}(t)\rho(t) ~~(x\neq 0)~~~~\mbox{and}~~~\hat f_0^{\{0\}}(t)=1-\rho(t)
 \eeqn
(cf.  \cite{U}).  We also have  $ q^k(0,y) =f^{\{0\}}_{-y}(k)$ as previously noted.  Now, employing
 the second identity in (\ref{111})  we see
\[
\sum_{k} r^{k}_\a(x,y) e^{ikt}=\frac{\ga\hat f^{\{0\}}_x(t)\hat f^{\{0\}}_{-y}(t)}{1+\ga \rho(t)}~~~~\mbox{with}~~~\ga=\frac{1-\a}{\a}.\]
 Hence, for all $x, y\in \Z$
$$r_\a^n(x,y)=\frac{\ga}{2\pi}\int_{-\pi}^\pi \frac{ \hat f^{\{0\}}_x(t)\hat f^{\{0\}}_{-y}(t)}{1+\ga \rho(t)}e^{-int}dt.$$

Since $\hat f^{\{0\}}_x(t)\hat f^{\{0\}}_{-y}(t)$ is the characteristic function of the convolution $f_x^{\{0\}}*f_{-y}^{\{0\}}$,
$$r_\a^n(x,y)=\ga f_x^{\{0\}}*f_{-y}^{\{0\}}(n)- \frac{\ga^2}{2\pi}\int_{-\pi}^\pi \frac{ \hat f^{\{0\}}_x(t)\hat f^{\{0\}}_{-y}(t)\rho(t)}{1+\ga \rho(t)}e^{-int}dt.$$
Under the constraints $|x|\vee|y|<a_\circ\sqrt n$ and  $\e\sqrt n< |x|\wedge |y| $ with $\e>0$  we apply Theorem A and,  making  scaling argument,  infer that  the first term  on the right side above agrees with that of the formula (\ref{rr}) with the factor $(|x|+|y|)\sigma^2 [a^*(x)a^*(-y)]/(a^*(x)+a^*(-y))|xy|$ multiplied, which
 factor  may be replaced by 1 since $|x|\wedge |y|\to \infty$. The second term is readily evaluated to be negligible  (see (\ref{113}) below and the argument following it as well as (\ref{114}),  which provides relevant  properties  of $\rho$).

By the first identity in (\ref{111}) we have 
\beqn\label{113}
\hat f^{\{0\}}_x(t)\hat f^{\{0\}}_{-y}(t)=\rho^2(t)\pi_{-x}(t)\pi_y(t)~~~~~\mbox{ if}~~~ x\neq 0, y\neq 0.
\eeqn
 As in Section 3 we write $\pi_{-x}(t)={\rm e}_x(t) +\pi_0(t)-a(x)$. Then \beq
\rho^2(t)\pi_{-x}(t)\pi_y(t)&=&\rho^2{\rm e}_{-x}{\rm e}_y+a(x)a(-y)\rho^2-(a(-y){\rm e}_x+a(x){\rm e}_{-y})\rho^2-1\\
&& +~ (\pi_{-x}+\pi_y)\rho.
\eeq
The contribution to $r_\a^n(x,y)$ of the last term  is $\ga(f_x^{\{0\}}(n)+f_{-y}^{\{0\}}(n))$. Those of the other terms are all
$o((|x|\vee|y|)/n^{3/2})$  if $|x|\wedge |y|=o(\sqrt n)$:
 for the proof we need the estimates  of ${\rm c}_x'''(t)$ and ${\rm s}_x'''(t)$ (which require no further moment condition) in addition to
those given in Lemmas B1 and B2. Now we may conclude that for $x\neq 0, y\neq 0$,
\beqn\label{112}
r_\a^n(x,y)=\frac{ 1-\a}{\a}(f_x^{\{0\}}(n)+f_{-y}^{\{0\}}(n))(1+o(1))~~~~~\mbox{as}~~~~\frac{|x|\wedge |y|}{\sqrt n}\to 0,
\eeqn
which agrees with (\ref{rr}) owing to  Theorem A again. If  $x= 0$, then $\hat f^{\{0\}}_x\hat f^{\{0\}}_{-y}=\rho(1-\rho)\pi_y=\rho\pi_y-\rho-\rho^2{\rm e}_{y}+\rho^2a(-y)$ and we readily obtain (\ref{112}). The case $y=0$ is similar.  ~~~\qed

\section{Appendix}
 ~ Under the basic assumption of this paper
Hoeffding \cite{H} shows that $\Re\{(1-e^{ixl})\phi^k(l)\}$  is summable on $\{k=1,2,\ldots\}\times \{-\pi<l<\pi\}$ and hence the series that defines $a(x)$ is absolutely convergent; and, as a consequence of it, that
\beqn\label{a(n)}
a(x)=  \frac1{2\pi} \int_{-\pi}^\pi \Re\bigg\{ \frac{1-e^{ixl}}{1-\phi(l)}\bigg\}dl
\eeqn 
(see \cite{K} and the references contained in  it for more information on $a(x)$).
 In this appendix we  derive an asymptotic estimate of $a(n)$  under the moment condition $E|Y|^{2+\de}<\infty$ ($0\leq \de\leq 2$). We also include a proof of (\ref{a(n)}) for reader's convenience.

Put for $0<|l|\leq \pi$,
$$ \phi_c(l)=\Re \phi (l)=E[\cos lY]~~~\mbox{ and}~~~~ \phi_s(l)=\Im \phi(l)=E[\sin lY].$$
Then
\beqn\label{a-1}
\Re\Bigg\{ \frac{1-e^{ixl}}{1-\phi(l)} \Bigg\}=\frac{1-\phi_c(l)}{|1-\phi(l)|^2} (1-\cos xl)+ \frac{\phi_s(l)}{|1-\phi(l)|^2} \sin xl.
\eeqn
Noticing $\phi_s(l)=E[\sin Yl -Yl]$, one infers that 
\beqn\label{s}
\int_{-\pi}^\pi \frac{|\phi_s(l)l|}{|1-\phi(l)|^2} \leq  E\int_{-\pi}^\pi \frac{|\sin Yl -Yl|}{|1-\phi(l)|^2}|l| dl\leq CE\bigg[Y^2\int_{0}^{\pi|Y|}\frac{|\sin u -u|}{u^3}du\bigg]<\infty;
\eeqn
hence by the dominated convergence theorem
\beqn\label{a-4}
\lim_{|x|\to\infty}\frac1{x}\int_{-\pi}^\pi     \frac{|\phi_s(l)|}{|1-\phi(l)|^2} |\sin xl |dl =0.
\eeqn
From this together with the equality  $\int_{-\infty}^\infty (1-\cos u)u^{-2}du =\pi$ 
we conclude the following result.

\begin{lem}\label{lem_a-1}~As $|x|\to\infty$
$$ \int_{-\pi}^\pi \Bigg|\Re\Bigg\{ \frac{1-e^{ixl}}{1-\phi(l)} \Bigg\}\Bigg|dl =O(x)
~~~\mbox{and} 
~~~ \frac1{2\pi}\int_{-\pi}^\pi \Re\Bigg\{ \frac{1-e^{ixl}}{1-\phi(l)}\Bigg\}dl = \frac{|x|}{\sigma^2}+o(x).$$
\end{lem}
\vskip2mm\n

The next proposition in particular implies the identity (\ref{a(n)}).
\begin{prop}\label{thm_a-2}~ With a  uniformly bounded term  $o_b(1)$ that tends to zero as $K\to\infty$
$$\sum_{k=0}^K\Big[p^k(0)-p^k(-x)\Big] = \frac1{2\pi} \int_{-\pi}^\pi \Re\bigg\{ \frac{1-e^{ixl}}{1-\phi(l)}\bigg\}dl (1+o_b(1)).$$
\end{prop}

\vskip2mm

{\it Proof.} ~  Since we have the  Tauberian  condition   $[p^k(0)-p^k(-x)]=o(1/k)$ (with $x$ fixed) as is assured by Lemma \ref{lem2.7}, owing to  the corresponding  Tauberian theorem  it suffices  (apart from the boundedness of convergence) to show  that the Abelian sum
 \beqn\label{AS}
 \frac1{x}\sum_{k=0}^\infty [p^k(0)-p^k(-x)]r^k  = \frac1{2\pi x}\int_{-\pi}^\pi\Re\Bigg\{\frac{1-e^{ixl}}{1-r\phi(l)}\Bigg\}dl
 \eeqn
converges as $r\uparrow 1$ to $1/x$ times the right side of (\ref{a(n)}).  The proof of this convergence  however is routine  and omitted.  From (\ref{46}) it is  clear  that  the convergence above is bounded uniformly  in $x$, which combined with the bound $|p^k(0)-p^k(-x)|\leq  C|x/k|$ (cf. Lemma \ref{lem2.7}) implies that the error term $o_b(1)$ is uniformly bounded (see  the proof of Tauber's theorem given in $\S$1.23 of \cite{T}). ~~ \qed
\v2\n
{\sc Remark.}~ 
The remainder term $o_b(1)$  in Proposition \ref {thm_a-2} does not uniformly (in $x$) approach zero  since 
the  Abelian sum in (\ref{AS})  tends to zero as $x\to\infty$ for each $r<1$. 

\v2\v2

If $E|Y|^3<\infty$, put
\beqn\label{b&C}
\la_3=\frac1{3\sigma^2} E[Y^3]~~~~\mbox{and}~~~~C^*=\frac1{2\pi}\int_{-\pi}^\pi \bigg[\frac{\sigma^2}{1-\phi(l)}-\frac1{1-\cos l}\bigg]dl,
\eeqn
where the integral of the imaginary part is understood to  vanish  because of  skew symmetry).
Since $\int_0^\pi\big[(1-\cos l)^{-1}-2l^{-2}\big]dl=2/\pi$, the constant $C^*$ may alternatively be given by
\beqn\label{C_def2}
C^*=\frac{\sigma^2}{2\pi}\int_{-\pi}^\pi\bigg[\Re\bigg\{\frac1{1-\phi(l)}\bigg\}-\frac2{\sigma^2l^2}\bigg]dl -\frac2{\pi^2}.
\eeqn
The integral in (\ref{C_def2}) as well as  the  real part of the integral in (\ref{b&C})   is absolutely convergent
 (under $E[|Y|^3]<\infty$).  In fact,  in the expression
\beqn\label{1-psi}
\Re\bigg\{\frac1{1-\phi(l)}\bigg\}=\frac{-\phi_s^2}{|1-\phi|^2(1-\phi_c)}+\frac1{1-\phi_c}
\eeqn
the first term on the right side  is bounded and, since $\phi_c-1+\frac12\sigma^2 l^2=E[\cos Yl-1+\frac12(Yl)^2]$,
\beqn\label{1-psi2}\int_{-\pi}^\pi \bigg[\frac1{1-\phi_c(l)}-\frac2{\sigma^2l^2}\bigg]dl \leq CE\bigg[|Y|^3\int_{0}^{\pi|Y|} \frac{\cos u-1+\frac12 u^2 }{u^4}du\bigg]<\infty. 
\eeqn

We write $\sg t\,= t/|t|$ $(t\neq 0)$. Suppose that $E |Y|^{2+\de}<\infty$ for some $0\leq \de\leq 2$.  
\vskip2mm
\begin{prop}\label{cor6.1}~~ ~ If $0\leq\de<1$, then $\sigma^2a(x)=|x|+o(|x|^{1-\de})$ (as $|x|\to\infty$) where  the error term is bounded if and only if $E[|Y|^3]<\infty$.
If $1\leq\de<2$, then  
$$\sigma^2a(x)=|x|+C^*-(\sg x)\la_3+o(|x|^{1-\de}).$$
If $\de=2$, this formula is valid with $o(|x|^{1-\de})$ replaced by $O(1/x)$.
\end{prop}
\vskip2mm\n
{\it Proof.}~  ~ In the case $\de=0$ the assertion is proved  both in \cite{H} and in \cite{S}, and  in fact immediate from (\ref{a(n)}) and Lemma \ref{lem_a-1}. The integral in  (\ref{a(n)}) may be written as 
\beqn\label{a000}
\sigma^2a(x)=\frac{\sigma^2}{2\pi}\int_{-\pi}^\pi\bigg[\frac1{1-\phi(l)}-\frac2{\sigma^2l^2}\bigg](1-e^{ixl})dl+
\frac{1}{2\pi}\int_{-\pi}^\pi\frac2{l^2}(1-e^{ixl})dl.
\eeqn
 The second term on the right side of  (\ref{a000})   equals
\beqn\label{2ndterm}
\frac{4}{2\pi}\bigg[|x|\int_{0}^\infty \frac{1-\cos l}{l^2}dl-\int_{\pi}^\infty \frac{1-\cos xl}{l^2}dl\bigg]=|x|-\frac2{\pi^2}+O(1/x^2).
\eeqn

Let $0<\de<1.$   The first term is then $o(|x|^{1-\de})$. In fact,  since $\phi(l)-1+\frac12 \sigma^2 l^2=o(|l|^{2+\de})$,
 for any $\e>0$ there exists a constant $M$ such that 
 $$\int_0^\pi\frac{|\phi(l)-1+\frac12 \sigma^2l^2|}{|1-\phi(l)|l^2}|1-e^{ixl}|dl\leq \e |x|^{1-\de}\int_0^\infty u^{-2+\de}|1-e^{iu}|du+M.$$ 
The  assertion concerning boundedness follows from Corollary \ref{lem2.5}.

Let $1\leq \de <2$. 
Then by (\ref{a000}), (\ref{C_def2}) and (\ref{2ndterm})  
$$\sigma^2 a(x)=|x|+C^*+I_c+I_s+O(x^{-2}),$$
where
$$I_c=-\frac{\sigma^2}{2\pi}\int_{-\pi}^\pi \bigg[\Re\bigg\{\frac1{1-\phi(l)}\bigg\}-\frac2{\sigma^2t^2}\bigg]\cos xl\, dl~~~~~~
\mbox{and}~~~~~~I_s=\frac{\sigma^2}{2\pi}\int_{-\pi}^\pi \frac{\phi_s}{|1-\phi|^2}\sin xl\, dl.$$
Recall $\phi_s(l)=  -\frac12\sigma^2\la_3l^3+o(|l|^{2+\de})$.  Then, employing a truncation argument  with  a smooth  cut-off function  $w(t)$ (i.e., $w$ vanishes outside $(-\pi,\pi)$ and equals  $1$ in a neighborhood of zero) along with  integration by parts one infers that
\beq
I_s&=&-\la_3\frac{\sigma^4}{2\pi}\int_{0}^\infty \frac{\sin xl}{\sigma^4l/4}dl+\int_{-\pi}^\pi {r(l)w(l)\sin xl}\,dl+o(1/x^3)\\
&=&-{\la_3\,{\rm sign}\,x}+o(|x|^{1-\de}).
\eeq
Here the remainder term $r(l)=o(|l|^{\de-2})$ with $ r'(l)=o(|l|^{\de-3})$.  Similarly we obtain $I_c=o(|x|^{1-\de})$ if $\de>1$.  If $\de=1$, an application of Riemann-Lebesgue lemma  with  (\ref{1-psi}) and (\ref{1-psi2}) taken into account  shows
$I_c =o(1)$.  The case $\de=2$ is similar and omitted.  ~~ \qed
\v2

From Subsection {\bf 2.1} we extract  the following result. 
\begin{prop} \label{cor2.2} ~~ Suppose that the walk is not left continuous. Then both $f_+(x)-x$ and $\sigma^2a(x)-f_+(x)$ are positive for all $x>0$ and tend to extended positive numbers as $x\to\infty$, which are finite if and only if $E[Y^3; Y>0]<\infty$. Moreover if  $E[Y^3; Y>0]=\infty$, then
$$\lim_{x\to\infty}\; \frac{f_+(x)-x}{\sigma^2a(x)- x}=\frac12.$$
\end{prop}
\vskip2mm\n
{\it Proof}.~  The assertions follow on combining Lemma \ref{lem2.50}   with (\ref{eq2.51}) and  (\ref{eq2.50}).  ~\qed
  
\begin{cor} \label{cor7.5} ~ Suppose that $E|Y|^3<\infty$. Then $\lim_{x\to \pm \infty} (\sigma^2 a(x)-|x|)= C^*\mp \la_3\geq 0\,;$ in particular $C^*=\la_3$ (resp. $-\la_3$)  if and only if the walk is left (resp. right) continuous.
\end{cor}
\vskip2mm\n
{\it Proof}.~ The first half follows from Proposition \ref{cor6.1} and the second one from the proposition above.  ~\qed
\v2\n
 
\vskip4mm\n
{\bf Acknowledgments.}~ The author wishes to thank the anonymous  referee for his carefully reading  the original manuscript and pointing out several significant errors therein.

  \v2\v2\n
  {\bf Note added in proof.}~
    Some asymptotic estimates of the transition probability of the walk killed on a half line are obtained in \cite{BD}, \cite{Car}, \cite{D},  \cite{VW} and \cite{A-D}. In the first four  papers  the problem is considered for  wider classes of random walks: the variance may be infinite with  the law of $Y$  in the domain of   attraction of a stable law  \cite{VW}, \cite{D} or  of the normal law \cite{BD} \cite{Car}; $Y$ is not necessarily restricted to the arithmetic variables. The very recent paper \cite{D}  describes  the asymptotic behavior of $p^n_{(-\infty,0]}(x,y)$ valid uniformly within the   region of stable deviation of the  space-time variables and it  in particular  contains    Theorem \ref{thm1.3} and Corollary \ref{cor1.1}  as a special case.  The result  for the region $a_\circ^{-1}\sqrt n<x, y <a_\circ \sqrt n$  is  contained also  in  \cite{BD}.
  Theorem \ref{thm1.3} with $x=1$ and  $y=o(\sqrt n)$  is  also a special case of Proposition 1 of \cite{BD}, Theorem 5 of  \cite{VW} and  readily derived from Theorem 2 of \cite{Car}. Similar results are proved  in  \cite{A-D}, where  Corollary \ref{cor1.1} is also obtained.   The methods used in these papers  fully rests on the Wiener-Hopf factorization or its extension recently developed and ours are quite different from them.

\end{document}